\documentclass[11pt]{article}

\usepackage{amssymb}
\usepackage{setspace} 
\usepackage{amsmath}
\usepackage{amsthm}
\usepackage{enumitem} 
\usepackage{graphicx}
\usepackage{tkz-euclide}
\usepackage{pgfplots}
\usepackage{physics}
\usepackage{amsfonts}
\usepackage{mathrsfs}
\usepackage{wrapfig}
\usepackage{subcaption,color,tikz-cd}

\usetikzlibrary{decorations.pathmorphing}
\usetikzlibrary{decorations.pathreplacing}

\newcommand{\Xt}{\tilde{X}}
\newcommand{\Proj}{\text{Proj}}
\newcommand{\gammat}{\tilde{\gamma}}
\newcommand{\Sigmat}{\tilde{\Sigma}}
\newcommand{\Sigmah}{\widehat{\Sigma}}
\newcommand{\taut}{\tilde{\tau}}
\newcommand{\gt}{\tilde{g}}
\newcommand{\etat}{\tilde{\eta}}
\newcommand{\psit}{\tilde{\psi}}
\newcommand{\Int}{\text{Int}}
\newcommand{\dom}{\text{dom}}
\newcommand{\range}{\text{range}}
\newcommand{\diam}{\text{diam}}
\newcommand{\cB}{\mathcal{B}}
\newcommand{\cD}{\mathcal{D}}

\newcommand{\cR}{\mathcal{R}}
\newcommand{\cS}{\mathcal{S}}
\newcommand{\V}{\mathcal{V}}
\newcommand{\Vt}{\tilde{\mathcal{V}}}
\newcommand{\A}{\mathcal{A}}
\newcommand{\At}{\tilde{\mathcal{A}}}

\usepackage{titlesec}
\numberwithin{equation}{section}

\usepackage[bookmarksnumbered=true]{hyperref} 
\hypersetup{
    colorlinks = true,
    linkcolor = black,
    anchorcolor = blue,
    citecolor = black,
    filecolor = blue,
    urlcolor = blue
    }

\newcommand{\Z}{\mathbb{Z}}
\newcommand{\R}{\mathbb{R}}

\newcommand\restr[2]{{
		\left.\kern-\nulldelimiterspace 
		#1 
		\vphantom{\big|} 
		\right|_{#2} 
}}

\titleformat{\section}
  {\normalfont\fontsize{15}{15}\bfseries}{\thesection}{1em}{}

\theoremstyle{definition}

\newtheorem{definition}{Definition}[section] 

\newtheorem{theorem}[definition]{Theorem}
\newtheorem{lemma}[definition]{Lemma}
\newtheorem{remark}[definition]{Remark}
\newtheorem{sublemma}[definition]{Sublemma}
\newtheorem{assumption}[definition]{Assumption}

\newtheorem{corollary}[definition]{Corollary}
\newtheorem{proposition}[definition]{Proposition}

\theoremstyle{plain} 
\newcommand{\thistheoremname}{}
\newtheorem*{genericthm}{\thistheoremname}
\newenvironment{namedthm}[1]
  {\renewcommand{\thistheoremname}{#1}%
   \begin{genericthm}}
  {\end{genericthm}}

\newcommand{\Addresses}{{
\bigskip
  \footnotesize

   \textsc{Mathematics and Computer Science Department, Wesleyan University, Middletown, CT}
	\par\nopagebreak
  \textit{E-mail address:} 
  \texttt{dconstantine@wesleyan.edu}
  \bigskip
  \footnotesize

   \textsc{Mathematics and Computer Science Department, Wesleyan University, Middletown, CT}
	\par\nopagebreak
  \textit{E-mail address:} 
  \texttt{eshrestha@wesleyan.edu}
  \bigskip
  \footnotesize

        \textsc{Department of Mathematics, Rice University, Houston, TX}
	\par\nopagebreak
  \textit{E-mail address:} 
  \texttt{yandi.wu@rice.edu}
  \bigskip
  \footnotesize
}}

\emergencystretch=\maxdimen
\begin{document}

\title{Sub-actions for geodesic flows on locally CAT(-1) spaces}
\author{David Constantine \and Elvin Shrestha \and Yandi Wu}
\date{}
\maketitle
\begin{abstract}
    We extend a result of Lopes and Thieullen on sub-actions for smooth Anosov flows to the setting of geodesic flow on locally CAT(-1) spaces. This allows us to use arguments originally due to Croke and Dairbekov to prove a volume rigidity theorem for some interesting locally CAT(-1) spaces, including quotients of Fuchsian buildings and surface amalgams.
\end{abstract}
\tableofcontents
\section{Introduction}\label{sec:intro}

The classical Liv\v{s}ic Theorem for a smooth, transitive Anosov flow $\{\phi_t\}$ states the following (\cite{livsic}): any H\"older function $A$ that integrates to zero over each closed orbit for the flow is itself a derivative. That is, there is a function $V$, smooth in the flow direction and still H\"older, such that $A(x)=\frac{d}{dt}|_{t=0}V(\phi_t x)$. (See, e.g., \cite{KH95} for a proof and discussion.)
Equivalently, but also providing a formulation for this statement which does not require $V$ to be smooth in the flow direction, $\int_0^T A(\phi_t x)dx = V(\phi^T x)- V(x).$

This theorem has far-reaching consequences for these flows. One with a particular connection to the present paper is that if two negatively curved Riemannian metrics on the same compact manifold have the same marked length spectrum, then their geodesic flows are conjugate. This is the starting point for Croke's proof of marked length spectrum rigidity for surfaces in \cite{croke}. This rigidity result was proved independently by Otal in \cite{otal} using tools which will play an important role in the current paper.

A natural generalization of the Liv\v{s}ic Theorem asks whether assuming that the periodic integrals of $A$ are all non-negative (or, with trivial modifications, non-positive) guarantees a $V$ whose derivative bounds $A$ below. Lopes and Thieullen term such a $V$ a `sub-action' for $A$ and prove that whenever $A$ is H\"older, a sub-action which is smooth in the flow direction and still H\"older exists \cite{LT}. Independently and concurrently Pollicott and Sharp proved a similar theorem in \cite{ps}. The proof in \cite{ps} is simpler, but establishes less -- they do not obtain that the sub-action $V$ is smooth in the flow direction or some of the more detailed results on regularity provided by \cite{LT}. See \cite{ps} for a good survey of related results in this area.

%
\subsection{Statement of results}\label{subsec:results}

In this paper, we follow the approach of \cite{LT}, applied specifically to geodesic flow on a CAT(-1) space. This is a \emph{metric Anosov flow} (or \emph{Smale flow}) -- it satisfies the essential properties of an Anosov flow, abstracted from the smooth setting to the general metric setting by Pollicott \cite{Pol87} (see \S\ref{subsec:metric anosov flows} below). 
This was proved in \cite{CLT20} where the authors found a coding for the geodesic flow using carefully chosen Poincar\'e sections. These sections turn out to be perfect candidates for the the sections used in the arguments of \cite{LT} (see \S\ref{sec:sections} below). Our main theorem is the following:

\begin{namedthm}{Main Theorem}\label{thm:flow subactions}
Let $(X, d_X)$ be a locally CAT(-1) space with geodesic flow $\{g_t\}$. Let $A: GX \to \R$ be  H\"older. Then there exists a map $V: GX \to \R$, called a sub-action, that is H\"older, smooth in the flow direction, such that for any geodesic $\gamma \in GX$ and every $T>0$, 
$$\int_0^T A \circ g_t(\gamma) \; dt \geq V \circ g_T(\gamma) - V(\gamma) + m T$$
for some constant $m = m(A)$. Equivalently, for any $\gamma \in GX $
$$A(\gamma) = m + \left( \frac{d}{dt}\right)\Big|_{t=0} V(g_t(\gamma)) + H(\gamma)$$
for some non-negative function $H: GX \to \R_{\geq 0}$ that is smooth along the flow direction and H\"older. 
\end{namedthm}

We turn to an application of the     \hyperref[thm:flow subactions]{Main Theorem} to the marked length spectrum in Section \ref{sec:MLS}. We prove a \textit{volume rigidity} result for \textit{surface amalgams}, which, roughly speaking, are constructed by identifying finitely many compact surfaces with boundary along their boundary components. For a more precise definition, we refer the reader to Section \ref{sec:surfaceamalgams}. Before stating the volume rigidity result for surface amalgams, we introduce some terminology.

\begin{definition}[Marked length spectrum] The \textit{marked length spectrum of a metric space $(X, g)$} is the class function 
$$\mathscr{L}_g: \pi_1(X) \rightarrow \mathbb{R}^{+}, [\alpha] \mapsto \inf\limits_{{\gamma \in [\alpha]}} \ell_g(\gamma)$$ 
which assigns to each free homotopy class $[\alpha] \in \pi_1(X)$ the infimum of lengths in the $g$-metric of curves in the class.
\end{definition}

Note that in the case of CAT($-1$) spaces, the marked length spectrum is simply a length assignment to every closed geodesic in $(X, g)$, as each homotopy class has a unique geodesic representative. 

Otal \cite{otal} and Croke \cite{croke} proved that for compact, negatively curved surfaces, if $\mathscr{L}_{g_0} = \mathscr{L}_{g_1}$ then $g_0$ and $g_1$ are isometric, with the classical Liv\v{s}ic theorem playing a role in Croke's proof. Suppose instead that $\mathscr{L}_{g_0} \leq \mathscr{L}_{g_1}$, i.e., for every free homotopy class $[\alpha] \in \pi_1(X)$, $\inf\limits_{\gamma \in [\alpha]} \ell_{g_0}(\gamma) \leq \inf\limits_{\gamma' \in [\alpha]} \ell_{g_1}(\gamma')$. In the setting of negatively curved surfaces, Croke and Dairbekov proved that such a marked length spectrum inequality implies a corresponding inequality in volumes of the surface with respect to $g_0$ and $g_1$. Furthermore, if the volumes are equal, then $g_0$ and $g_1$ are isometric \cite{CD04}. The proof of rigidity in the equality case crucially uses Lopes and Thieullen's sub-action result. Having extended that result to the CAT(-1) setting, we use Croke and Dairbekov's idea as well as some recent results on marked length spectrum rigidity for surface amalgams and similar spaces to prove:

\begin{namedthm}{Volume Rigidity Corollary}\label{cor:volumerigiditySA}
Let $(X, g_0)$ and $(X, g_1)$ be two simple, thick, negatively curved surface amalgams satisfying certain smoothness conditions around the gluing curves. Suppose $\mathscr{L}_{g_0} \leq \mathscr{L}_{g_1}$. Then Vol$_{g_0}(X) \leq \text{Vol}_{g_1}(X)$. Furthermore, if Vol$_{g_0}(X) = \text{Vol}_{g_1}(X)$, then $(X, g_0)$ and $(X, g_1)$ are isometric. 
\end{namedthm}

The fact that the sub-action $V$ is smooth in the flow direction is used in the proof of the \hyperref[cor:volumerigiditySA]{Volume Rigidity Corollary}. This provides some justification for our adoption of \cite{LT}'s more complicated but slightly stronger proof strategy.

%
\subsection{An outline of the paper}\label{subsec:outline}

In Section \ref{sec:preliminaries}, we collect definitions and basic results about CAT(-1) spaces, geodesic flow on such spaces, and its properties. We also prove some basic geometric facts which will be used later in the paper.

In Section \ref{sec:sections} we describe the construction of Poincar\'e sections for the geodesic flow which allow us to `discretize' the flow. We describe and prove some key properties of these sections that will be useful for subsequent arguments.

The proof of the \hyperref[thm:flow subactions]{Main Theorem} begins in earnest with Section \ref{sec:discretized}. In this section, we follow the arguments of \cite{LT} and solve the discretized version of the sub-action problem provided by the Poincar\'e sections. 

Section \ref{sec:flow} continues in the steps of \cite{LT} to extend the solution of the discretized problem to a sub-action for the flow. A careful inductive scheme allows one to make this extension while ensuring the desired regularity of the sub-action.

Finally, in Section \ref{sec:MLS} we use the \hyperref[thm:flow subactions]{Main Theorem} to prove the \hyperref[cor:volumerigiditySA]{Volume Rigidity Corollary} in its full generality.\\

\noindent \textbf{Acknowledgements.} The third author gratefully acknowledges generous support from an NSF RTG grant, DMS-2230900, during her time at University of Wisconsin, Madison. The authors would also like to thank Nicola Cavallucci for input on doubling metric spaces, which aided the development of the ideas in Section 2.4.
\section{Preliminaries}\label{sec:preliminaries}

%

\subsection{Geodesics and Geodesic Flow}\label{subsec:geos and geo flow}

Let $(\Xt, d_{\Xt})$ be a CAT(-1) space and $\Gamma$ be a discrete group of isometries of $\Xt$ acting freely, properly discontinuously, and cocompactly. (See \cite{bh} for background on CAT(-1) spaces.) The resulting quotient $X=\Xt/\Gamma$, with the metric $d_X$ induced by $d_{\Xt}$, is a compact, locally CAT(-1) space.

To study the geodesic flow on $X$, we need an analogue of the unit tangent bundle appropriate for this non-smooth setting which admits a unit-speed geodesic flow. This is provided by the following definition.

\begin{definition}\label{defn:geo flow}
    The space of (unit-speed) geodesics of $\Xt$ is
    $$G\Xt := \{\gammat: \R \to \Xt: d_{\Xt}(\gammat(s), \gammat(t)) = |s-t|\}.$$
    The geodesic flow $g_t$ on $G\Xt$ is given by $g_t \gammat(s) = \gammat(s+t)$ for any $t \in \R$. 
    $GX = G\Xt / \Gamma$ is the space of geodesics of $X$.
\end{definition}

There are natural metrics on these spaces.

\begin{definition}\label{defn:GX metric}
We equip $G\Xt$ with the metric 
$$d_{G\Xt}(\gammat_1, \gammat_2) := \int_{-\infty}^{\infty} d_{\Xt}(\gammat_1(t), \gammat_2(t)) e^{-2|t|} \; dt,$$
and $GX$ with the metric
$$d_{GX}(\gamma_1, \gamma_2) = \min\limits_{\widetilde{\gamma_1}, \widetilde{\gamma_2}}\int_{-\infty}^{\infty} d_{\widetilde{X}}(\widetilde{\gamma_1}(t) , \widetilde{\gamma_2}(t))e^{-2\lvert t \rvert}dt$$ 
\noindent where the minimum is taken over all lifts $\widetilde{\gamma_1}$ and $\widetilde{\gamma_2}$ of $\gamma_1$ and $\gamma_2$, respectively. 
\end{definition}

A straightforward computation shows that the geodesic flow is unit-speed with respect to $d_{GX}$. (This is the reason for the normalizing factor 2 in the exponent.) We record here a few basic facts about $d_{GX}$ and its interaction with the geodesic flow.

\begin{lemma}\label{lem:dGX to dX}\cite[Lemma 2.8]{clt}
For all $\gamma_1, \gamma_2 \in GX$, $d_X(\gamma_1(0), \gamma_2(0))\leq 2d_{GX}(\gamma_1, \gamma_2)$.
\end{lemma}

\begin{lemma}\label{lem:geo flow is Lipschitz}\cite[Lemma 2.5]{CLT20}
For any $t$ and any $\gamma_1,\gamma_2\in GX$, $d_{GX}(g_t\gamma_1, g_t\gamma_2) \leq e^{2|t|}d_{GX}(x,y)$.
\end{lemma}

A central tool in the geometry of CAT(-1) spaces is the boundary at infinity.

\begin{definition}\label{defn:boundary at infinity}
The\textit{ boundary at infinity }$\partial^\infty \Xt$ is the space of equivalence classes of geodesic rays in $\Xt$, where two rays $c_1, c_2: [0, \infty) \to \Xt$ are equivalent if there exists $M \geq 0$ such that $d_{\Xt}(c_1(t), c_2(t)) \leq M$ for all $t\geq 0$. For a geodesic $\gammat$ we denote by $\gammat(-\infty)$ and $\gammat(+\infty)$ its backward and forward endpoints at $\partial^\infty \Xt$. 
\end{definition}

In CAT(-1) spaces, a pair $\{\xi, \eta\}$ of distinct points on $\partial^\infty \Xt$ uniquely determines an unparameterized, unoriented geodesic. Once we specify an orientation and the time-$0$ point on the geodesic, we have a geodesic in $G\Xt$. Therefore $G\Xt$ can be identified with $[(\partial^\infty \Xt \times \partial^\infty \Xt) \backslash \Delta] \times \R$, where $\Delta$ is the diagonal. For ease of notation, let $\partial^{(2)}_\infty\Xt = (\partial^\infty \Xt \times \partial^\infty \Xt) \backslash \Delta.$

%

\subsection{Stable and Unstable Sets}\label{subsec:stable unstable}

Under the geodesic flow on a negatively curved manifold, the unit tangent bundle admits foliations by stable and unstable manifolds which are essential tools for studying the dynamics of these flows. In the CAT(-1) setting, we have a purely geometric description of analogous stable and unstable sets. Below, we will note that these are also analogous in their dynamical role.

\begin{definition}\label{defn:busemann}
    Let $p \in \Xt$, $\xi\in \partial^\infty\Xt$, and $\gammat$ be the geodesic  ray from $p$ to $\xi$. The \textit{Busemann function centered at $\xi$ with basepoint $p$ }is defined as
    \begin{align*}
        B_p(-, \xi): \Xt &\to \R \\
        q & \mapsto \lim_{t\to \infty} \left(d_{\Xt}(q, \gammat(t)) - t\right).
    \end{align*}
    For convenience, given geodesic ray $\gammat(t)$, the \textit{Busemann function determined by} $\gammat$ is the Busemann function centered at $\gammat(+\infty)$ with basepoint $\gammat(0)$: 
    $$B_{\gammat}(-):= B_{\gammat(0)}(-, \gammat(+\infty)).$$
    The level sets for $B_p(-, \xi)$ are \textit{horospheres}.
\end{definition}

\begin{definition}\label{def:stablesetsBusemann}
Let $\gammat \in G\Xt$.
\begin{itemize}
    \item The \textit{strong stable set        through} $\gammat$ is 
        $$W^{ss}(\gammat) = \{\gammat'\in G\Xt: \gammat'(+\infty) = \gammat(+\infty) \text{ and } B_{\gammat}(\gammat'(0)) = 0 \}.$$

    \item     For $\delta>0$,
        \begin{align}
            W^{ss}_\delta(\gammat) = \{\gammat'\in G\Xt : \ & \gammat'(+\infty) = \gammat(+\infty),\nonumber \\
                & B_{\gammat}(\gammat'(0)) = 0, \text{ and } d_{G\Xt}(\gammat, \gammat') < \delta\}.\nonumber
        \end{align}

    \item The \textit{weak stable set through }$\gammat$ is 
        $$W^{cs}(\gammat) = \bigcup_{t \in \R} g_t(W^{ss}(\gammat)).$$

    \item Similarly, the \textit{strong and weak unstable sets through} $\gammat$ are 
    $$W^{uu}(\gammat) = \{\gammat'\in G\Xt: \gammat'(-\infty) = \gammat(-\infty) \text{ and } B_{-\gammat}(\gammat'(0)) = 0 \},$$
    \begin{align}
        W^{uu}_\delta(\gammat) = \{\gammat'\in G\Xt : \ & \gammat'(+\infty) = \gammat(+\infty),\nonumber \\
            & B_{-\gammat}(\gammat'(0)) = 0 \text{ and } d_{G\Xt}(\gammat, \gammat') < \delta\}, \nonumber
    \end{align}
    and
    $$W^{cu}(\gammat) = \bigcup_{t \in \R} g_t(W^{uu}(\gammat)).$$
\end{itemize}
\end{definition}

\begin{figure}[h]
\centering
\includegraphics[width=0.6\textwidth]{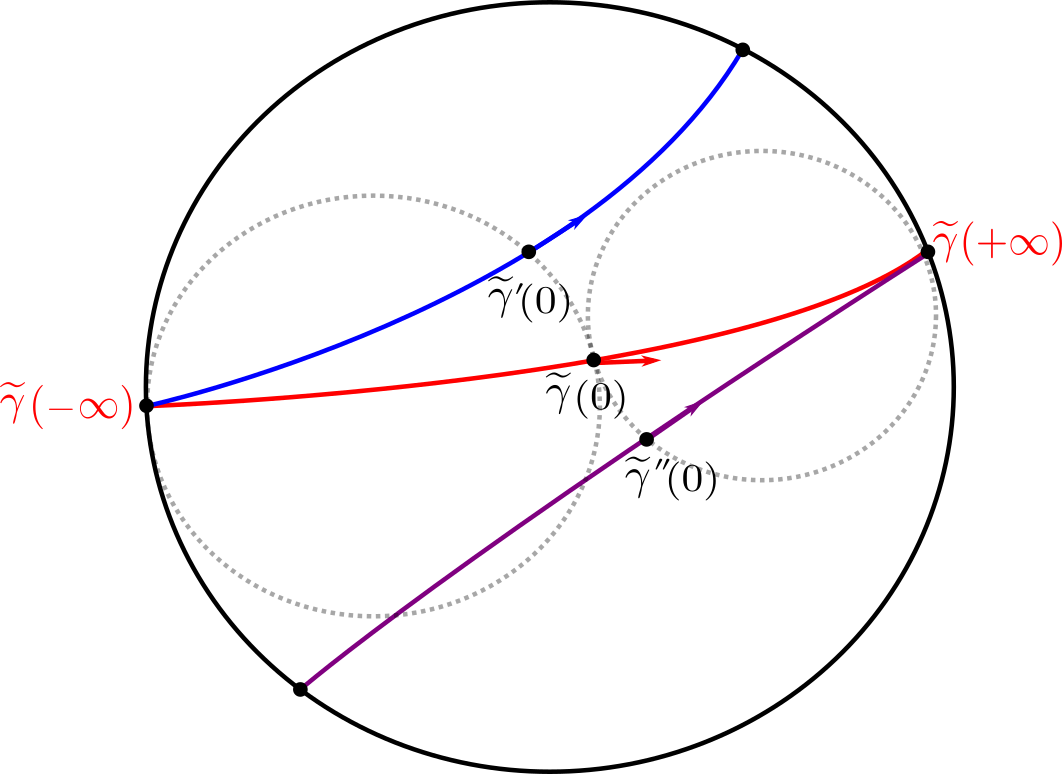}
\caption{In the above picture, $\widetilde{\gamma'} \in W^{uu}(\widetilde{\gamma})$ and $\widetilde{\gamma''} \in W^{ss}(\widetilde{\gamma})$. The 0-level sets of the Busemann functions are shown in gray.}
\end{figure}

The next lemma outlines some useful, standard properties of Busemann functions. Their proofs are straightforward.

\begin{lemma}\label{lem:BusemannProperties}
    Let $\gammat \in G\Xt$.
    \begin{itemize}
        \item[(1)] For any $s \in \R$, $B_{g_s\gammat}(-) = B_{\gammat}(-) + s$, so the $0$-level set for $B_{g_s\gammat}(-)$ is the $(-s)$-level set for  $B_{\gammat}(-)$. 
        \item[(2)] If $\etat \in W^{cs}(\gammat)$, then for $s_1, s_2 \in \R$, $B_{\gammat}(\etat(s_1))- B_{\gammat}(\etat(s_2)) = s_2 - s_1$.
        \item[(3)] $B_{\gammat}(-): \Xt \to \R$ is $1$-Lipschitz. 
    \end{itemize}
\end{lemma}

%

\subsection{Metric Anosov Flows}\label{subsec:metric anosov flows}

In \cite{Pol87}, Pollicott defines metric Anosov flows, generalizing the essential properties of Anosov flows to the non-smooth setting.

A continuous flow on a compact metric space is Anosov if, roughly speaking, it has a topological local product structure which coheres with the dynamics of the flow in a way which mimics the analogous structure for Anosov flows. Locally, $Y$ looks like a product -- two nearby points can always be connected by a small step along a `stable' set, then a small step along an `unstable' set, then a small move along the flow. 
This coheres with the dynamics of the flow in that the `stable' set is truly stable in the dynamical sense: there are constants $C,\lambda>0$ such that if $x$ and $y$ are nearby points on the same stable set, then
\[d(\phi_t x, \phi_t y) \leq Ce^{-\lambda t}d(x,y) \ \mbox{ for } t\geq 0.\]
The corresponding expression for exponential contraction in backwards time holds for pairs of nearby points on a common unstable set. (See \cite{Pol87} or \cite{CLT20} for more detail.)

Metric Anosov flows have many of the properties of Anosov flows. As a first example (which we will use below), we have:

\begin{proposition}[\cite{Bow73}, Cor 1.6 and \cite{Pol87}, Prop 1]\label{theorem:anosovflowexpansive}
A metric Anosov flow is expansive.
\end{proposition}

The flows we consider in this paper are metric Anosov:

\begin{theorem}[\cite{CLT20}, Theorem 3.4]\label{theorem:GeodesicFlowMetricAnosov}
For a compact, locally CAT(-1) space $X$, the geodesic flow on $GX$ is metric Anosov with $\lambda=1$. Specifically, the strong stable and unstable sets of Definition \ref{def:stablesetsBusemann} are the stable and unstable sets for the Anosov flow. That is, there exist small $\delta>0$ and $C>0$ such that
\[d_{G\Xt}(g_t \gammat, g_t \gammat') \leq Ce^{-t} d_{G\Xt}(\gammat,  \gammat') \mbox{ for } t\geq 0 \mbox{ and }\gammat' \in W^{ss}_\delta(\gammat);\]
\[d_{G\Xt}(g_{-t} \gammat, g_{-t} \gammat') \leq Ce^{-t} d_{G\Xt}(\gammat,  \gammat') \mbox{ for } t\geq 0 \mbox{ and }\gammat' \in W^{uu}_\delta(\gammat).\]
\end{theorem}

%

\subsection{Doubling and covering properties}\label{subsec:doubling and covering}

In the proof of Theorem \ref{thm:flow subactions}, we will need a technical fact about the geometry of $GX$. We collect the necessary arguments leading to that fact here.

Let $(X, d)$ be a metric space. We introduce two useful properties of metric spaces:

\begin{definition}\label{defn:doubling}
An open set $\mathcal{U} \subset X$ with compact closure satisfies the \textit{doubling} property if there exists $N_1 \in \mathbb{N}$ such that for all $r > 0$ and $x \in \mathcal{U}$, $\overline{B(x, 2r)}$ can be covered by at most $N_1$ balls of radius $r$. 
\end{definition}

\begin{definition}\label{defn:packing}
An open set $\mathcal{U} \subset X$ with compact closure satisfies the \textit{packing} property if there exists $N_2 \in \mathbb{N}$ such that for all $r > 0$ and $x \in \mathcal{U}$, the maximum cardinality of an $r$-separated subset of $B(x, 2r)$ is less than or equal to $N_2$. 
\end{definition} 

One can prove that doubling implies packing (although the converse is not true).

\begin{lemma}\label{lem:doublingpacking}
If $\mathcal{U}$ satisfies the doubling property, then it also satisfies the packing property. 
\end{lemma}

\begin{proof}
Take an arbitrary $x \in \mathcal{U}$ and $r > 0$, and suppose $\overline{B(x, 2r)}$ can be covered by at most $N_1$ balls of radius $r$, and these balls can in turn be covered by $N_1$ balls of radius $\frac{r}{2}$. We thus have a covering $\{B(x_i, \frac{r}{2})\}_{i \in I}$ of $\mathcal{U}$ of cardinality at most $N_1^2$. Consider an $2r$-separated set $S(2r) = \{x_i \in \overline{B(x, 2r)} : d(x_i, x_j) \geq 2r \text{ for all $i \neq j$}\}$. Note that each $B(x_i, \frac{r}{2})$ can contain at most one element of $S(2r)$; otherwise, if $y_1, y_2 \in B(x_i, \frac{r}{2})$, then $d(y_1, y_2) \leq r$. Hence $\mathcal{U}$ satisfies the packing property with $N_2=N_1^2$.
\end{proof}

Now let $X$ be locally CAT(-1). Provided the double boundary $\partial_\infty^{(2)}\Xt$ of $\widetilde{X}$ satisfies the doubling property, one can show an analogue of a technical lemma (Lemma 22) in \cite{LT}.

\begin{lemma}
\label{lem:covering}
    Suppose that $\partial_\infty^{(2)}\Xt$ satisfies the doubling property for some metric. Then there exists a constant $C$ depending only on $X$ such that for every $\epsilon > 0$ and every open $\Sigma \subset \partial_\infty^{(2)}\Xt$ with $\overline{\Sigma}$ compact, there exists a finite cover $\{B_{jk}\}_{j\in J, 1\leq k\leq C}$ of $\Sigma$ by $\epsilon$-balls such that for every fixed $1 \leq k^{\ast} \leq C$, the balls in $\{B_{jk^{\ast}}\}_{j \in J}$ are disjoint. 
\end{lemma}

\begin{proof}
    Let $\{x_l\}_{l \in L}$ be a maximal set in $\Sigma$ such that $d(x_l, x_m) \geq 2 \epsilon$. By maximality, $\bigcup_{l \in L} B(x_l, 2\epsilon)$ covers $\Sigma$. By the doubling condition, we can cover each $B(x_l, 2\epsilon)$ by at most $N_1$ $\epsilon$-balls, for some $N_1$ depending only on $\Sigma$. For each $l \in L$, we choose this collection minimally and use $\{B(y_{lm}, \epsilon)\}_{m \leq M(l)}$ to denote the collection covering $B(x_l, 2\epsilon)$. 
    
    Note that if $B(y_{l_1m_1}, \epsilon) \cap B(y_{l_2m_2}, \epsilon) \neq \varnothing$, then $d(x_{l_1}, x_{l_2}) \leq 6\epsilon$. Indeed, if $x \in B(y_{l_1m_1}, \epsilon) \cap B(y_{l_2m_2}, \epsilon)$, then for $i = 1, 2$, we have that $d(y_{l_im_i}, x_{l_im_i}) \leq 2\epsilon$ and $d(x, y_{l_im_i}) < \epsilon$, so by the triangle inequality, $d(x_{l_1}, x_{l_2}) \leq 6\epsilon$. Furthermore, by Lemma \ref{lem:doublingpacking}, since $d(x_{l_1}, x_{l_2}) \geq 2\epsilon$, for each $x_{l_1}$, we have:    
    \begin{align}
    \#\{y_{l_2m_2}: B(y_{l_1m_1}, \epsilon) \cap B(y_{l_2m_2}, \epsilon) \neq \varnothing\} &\leq N_1\big(\#\{x_{l_2}: d(x_{l_1}, x_{l_2}) \leq 6\epsilon\}\big) \nonumber \\
    &\leq N_1N_2^2 = N_1^5. \label{N2}
    \end{align} 
    
    We now use $\{B(y_{l_m})\}_{l \in L}$ to construct a finite-valence dual graph $\Gamma$ as follows:
    
    \begin{itemize}
        \item $V(\Gamma) = \{y_{l_m}\}_{l \in L, 1 \leq m \leq M(l)}$ (where $M(l) \leq N_1$);
        \item $E(\Gamma) = \{[y_{l_1m_2}, y_{l_2m_2}] : B(y_{l_1m_1}, \epsilon) \cap B(y_{l_2m_2}, \epsilon) \neq \varnothing \}$. 
    \end{itemize}

By Equation \eqref{N2}, each vertex in $V(\Gamma)$ has valence at most $N_1N_2^2$. By Brooks' Theorem (\cite{brooks} or \cite{lovasz}), $\Gamma$ can be colored with at most $N_1N_2^2+1$ colors. Relabel the set $\{y_{lm}\}$ to $\{y_{jk}\}$ where $k \in [1, C]$ indexes the colors in the coloring of $V(\Gamma)$, and $j \in [1, J(k)]$ indexes the vertices colored with the $k$th color. 

By definition of a coloring of a graph, for a fixed $k^{\ast}$, the set $\{y_{jk^{\ast}}\}_{1 \leq j \leq J(k^{\ast})}$ is a totally disconnected subgraph of $\Gamma$. Then the balls in $\{B(y_{jk^{\ast}}, \epsilon)\}_{1 \leq j \leq J(k^{\ast})}$ are pairwise disjoint, as desired. 
\end{proof}

To use Lemma \ref{lem:covering}, we need a metric on $\partial_\infty^{(2)}\Xt$ which is doubling. Recall that the \textit{$\ell^{\infty}$ product metric} on the product $X \times Y$ of two metric spaces $(X, d_X)$ and $(Y, d_Y)$ is defined by:
$$d_{X \times Y}\big((x_1, y_1), (x_2, y_2)\big) = \max\{d_X(x_1, x_2), d_Y(y_1, y_2)\}.$$

The proof of the following is relatively straightforward. 

\begin{lemma}\label{product-doubling-lemma} The $\ell^{\infty}$ product metric space $(X \times Y, d_{X \times Y})$ is doubling if $(X, d_X)$ and $(Y, d_Y)$ are doubling. More precisely, if $M$ and $N$ are the doubling constants of $(X, d_X)$ and $(Y, d_Y)$ respectively, then $MN$ is the doubling constant of $(X \times Y, d_{X \times Y})$. 
\end{lemma}


Recall that a metric space is \textit{proper} if every closed ball is compact. A group $G$ acting on $X$ by isometries is \textit{discrete} if for every ball $B = B(x, r) \subset X$, $\{g \in G: gB \cap B \neq \varnothing\}$ is finite. We now show that for a proper geodesic Gromov hyperbolic metric space $X$ equipped with a discrete, cocompact action (e.g. the CAT($-1$) space $\Xt$), $\partial_\infty^{(2)}\Xt$ equipped with a ($\ell^{\infty}$) product of \textit{visual metrics} is doubling. In order to do so, we must define some stronger properties of metric spaces. 

\begin{definition} A metric space $X$ is \textit{Ahlfors Q-regular} if there exists a positive Borel measure $\mu$ on $X$ such that: $$C^{-1}R^Q \leq \mu\big(B(x, R)\big) \leq CR^Q.$$
\end{definition} 

While we defined doubling for metric spaces, we can also say a \textit{measure} admitted by a metric space is doubling. 

\begin{definition} A positive Borel measure $\mu$ on a metric space $X$ is \textit{doubling} if there exists a constant $C_{\mu}$ such that:
$$\mu(B(x, 2r)) \leq C_{\mu}\mu(B(x, r)).$$
\end{definition}

We now briefly prove a lemma that is usually assumed in the literature due to the simplicity of the proof: 

\begin{lemma} \label{lemma:ahlfors-doubling} Any Ahlfors Q-regular metric space satisfies the doubling property. 
\end{lemma}

\begin{proof}
We first note that if a metric space is Ahlfors Q-regular, then it admits a doubling measure. Indeed, $\mu$ itself is a doubling measure; by the Ahlfors $Q$-regularity of $X$, we have:
$$C^{-1}R^Q \leq \mu(B(x, R)) \leq CR^Q \text{ and } C^{-1}2^QR^Q \leq \mu(B(x, 2R)) \leq C2^QR^Q.$$
Thus, $R^Q \leq C\mu(B(x, R))$, so: 
$$\mu(B(x, 2R)) \leq C2^QR^Q \leq 2^QC^2\mu(B(x, R)).$$
Thus, $\mu$ is doubling with constant $2^QC^2$. Finally, by \cite{ljs}, we have that a metric space is doubling if and only if it admits a doubling measure. 
\end{proof} 

We equip the boundary of a CAT($-1$) space with a \textit{visual metric}. First, if $X$ is CAT($-1$), given a choice of basepoint $x_0 \in X$, we define the \textit{Gromov product} of $\xi, \eta \in \partial^{\infty} X$, which is denoted by $(\xi.\eta)_{x_0}$: 
$$(\xi.\eta)_{x_0} = \lim\limits_{\substack{x_n \rightarrow \xi \\ y_n \rightarrow \eta}} \frac{1}{2}\big(d(x_0, x_n) + d(x_0, y_n) - d(x_n, y_n)\big).$$

\begin{definition}[Visual metric for CAT(-1) spaces, \cite{bourdon}] If $X$ is a \textit{proper} CAT(-$1$) space, then we can define the \textit{visual metric} to be:
$$d_e(\eta, \xi) = \begin{cases}e^{-(\xi.\eta)_{x_0}} &\text{ if $\xi \neq \eta$} \\ 0 &\text{otherwise}.\end{cases}$$
\end{definition} 

Although we do not delve into the definition of visual metrics for Gromov hyperbolic spaces in general, we remark that they exist (see \cite{bh} Chapter III.H). This allows us to state the following theorem: 

\begin{theorem}[\cite{kleiner} Theorem 3.3]\label{ahlfors-thm} If $X$ is a proper, geodesic Gromov hyperbolic metric space which admits a discrete, cocompact isometric action, then $\partial^{\infty} X$ equipped with a visual metric is Ahlfors $Q$-regular for some $Q$.
\end{theorem}

The original proof of Theorem \ref{ahlfors-thm} is actually due to Coornaert (see \cite{coornaert}, Proposition 7.4), who shows that any $\Gamma$-quasiconformal measure on $\partial^{\infty} X$ with support contained in the limit set of $\Gamma$ is Ahlfors $Q$-regular, where $\Gamma$ is a discrete, cocompact isometric group action on $X$. In particular, if $\mu$ is any non-zero measure on $\partial^{\infty} X$, then $\mu$ is $\Gamma$-quasiconformal (see \cite{coornaert} Proposition 4.3). As a consequence of Lemma \ref{lemma:ahlfors-doubling}, we thus have: 

\begin{corollary}
If $X$ is a proper, geodesic Gromov hyperbolic metric space which admits a discrete, cocompact isometric action, then $\partial^{\infty} X$ equipped with a visual metric is doubling.     
\end{corollary}

Finally, we reach a key corollary: 

\begin{corollary}
If $X$ is a locally CAT($-1$) space, then $\partial^{(2)}_{\infty}(\widetilde{X})$ equipped with the $\ell^{\infty}$ product of two visual metrics is doubling.     
\end{corollary}

\begin{proof}
 Note that by properties of locally CAT($-1$) spaces, $\widetilde{X}$ is a proper geodesic Gromov hyperbolic space that admits a discrete, cocompact isometric group action, so by Theorem \ref{ahlfors-thm}, $\partial^{\infty} \widetilde{X}$ equipped with the visual metric is doubling. By Lemma \ref{product-doubling-lemma}, $(\partial^{\infty} \widetilde{X} \times \partial^{\infty} \widetilde{X}) \setminus \Delta = \partial^{(2)}_{\infty} \widetilde{X}$ equipped with the $\ell^{\infty}$ product of two visual metrics is doubling.
\end{proof}


\section{Poincar\'e sections and symbolic codings for the geodesic flow}\label{sec:sections}

A key tool in \cite{LT}'s approach to the sub-actions problem is to introduce carefully chosen Poincar\'e sections which discretize the geodesic flow, encoding it with a shift space. We follow the same approach, using the work already done in \cite{CLT20} to produce a particularly nice collection of sections which satisfies the main conditions needed in \cite{LT}'s argument.

We begin by recalling the constructions and properties of these sections from \cite{CLT20}. Then we check that the key properties needed in the main argument of this paper hold for these sections.

%

\subsection{Sections and Markov proper families}\label{subsec:markov families}

Again, let $(Y, d_Y)$ be a compact metric space with a continuous flow $\{\phi_t\}$. Assume that $\phi_t$ has no fixed points.

\begin{definition}\label{def:section}
A \textit{Poincar\'e section} (or simply \textit{section}) is a closed subset $D \subset Y$ such that the map $(y, t) \mapsto \phi_t y$ is a homeomorphism between $D \times [-\tau^*, \tau^*]$ and $\phi_{[-\tau^*, \tau^*]}D$ for some time $\tau^* > 0$.
\end{definition} 

We will use the following notation related to sections.

\begin{definition}\label{def:interior and projection}
Given a section $D$, we denote by $\Int_{\phi}D$ the interior of $D$ transverse to the flow; that is,
$$\Int_\phi D = D \cap \bigcap_{\epsilon>0}\left(\phi_{(-\epsilon, \epsilon)} D \right)^\circ.$$
Let $\Proj_D: \phi_{[-\tau, \tau]}D \to D$ be the projection map defined by $\Proj_D(\phi_t y) = y$.
\end{definition}

Given a collection $\mathcal{D} =\{D_1, \ldots, D_n \}$ of disjoint sections, if $\bigcup_{i=1}^n \phi_{(-\alpha,0)}\Int_\phi D_i =Y$, then any orbit of $\phi_t$ crosses an infinite sequence of sections. We let $\psi$ be the first-return map for this collection, and let $\tau$ be the first return time for the collection, defined by $\psi(x)=\phi_{\tau(x)}x$ for all $x\in \bigcup_{i=1}^n D_i$.

Any orbit of $\phi_t$ can be encoded by the bi-infinite sequence of sections it crosses. To make this encoding useful, however, sections need to be chosen carefully. Bowen (\cite{Bow73}) and Pollicott (\cite{Pol87}) provide a framework for doing this in the Anosov and metric Anosov flow settings, respectively. \cite{CLT20} further describes this process. The end result is a Markov proper family:

\begin{definition}\label{def:Markov proper family}
    Let $\cB = \{B_1, \ldots, B_n\}$ and $\cD = \{D_1, \ldots, D_n\}$ be collections of disjoint sections. $(\cB,\cD)$ form a \textit{Markov proper family at scale $\alpha>0$} if they satisfy the following properties:
\begin{itemize}
    \item[(1)] $\diam(D_i)<\alpha$ and $B_i \subset D_i$ for each $ i = 1, \ldots, n$;
    \item[(2)] $\bigcup_{i=1}^n \phi_{(-\alpha, 0)}\Int_{\phi} B_i = Y$;
    \item[(3)] Each $B_i$ is a \emph{rectangle} in $D_i$ (see \cite[Defn 3.6]{CLT20} for details);
    \item[(4)] $\cB$ satisfies a \emph{Markov property} (see \cite[Defn 3.8]{CLT20} for details).
\end{itemize}
\end{definition}

\cite{Bow73}, \cite{Pol87}, and \cite{CLT20} all discuss conditions under which Markov proper families exist. Condition (3) does not play a role in the present paper, and we will work with a weaker version of condition (4) (see Lemma \ref{lem:markov} below). For our purposes, the key result is the following.

\begin{proposition}[\cite{CLT20}, \S4]\label{prop:MarkovProperFamiliesExist}
Let $X$ be a compact, locally CAT(-1) space. For geodesic flow on $GX$, at any scale $\alpha>0$, there exist Markov proper families with Lipschitz return time functions $\tau$.
\end{proposition}



%
\subsection{Markov Codings and Markov Proper Families}

In this section, let $(\cR, \cS)$ be a Markov proper family at scale $\alpha>0$ for a metric Anosov flow $\{\phi_t\}$ on $Y$. We discuss how the sections index orbits as sequences in a subshift of finite type. 

\begin{definition}\label{def:CanonicalCodingSpace}
    For a Markov proper family $(\cR, \cS)$ for a metric Anosov flow, we define the \textit{canonical coding space} to be 
    $$\Theta = \Theta(\cR) = \left\{\theta \in \prod_{-\infty}^\infty\{1, \ldots, N\}: \text{ for all } l, k \geq 0, \bigcap_{j=-k}^{l} \psi^{-j}(\Int_\phi R_{\theta_j}) \not= \emptyset\right\}.$$
\end{definition}

In Section 2.3 of \cite{Pol87}, the symbolic space $\Theta(\cR)$ is shown to be a shift of finite type. 
Let $\sigma: \Theta \to \Theta$ be the associated left shift map.

\begin{definition}\label{def:stableunstablefortheta}
    For any $\theta \in \Theta$, we define the local stable and unstable sets of $\theta$ to be:
    \begin{align*}
        W_{loc}^s(\theta) & := \left\{y \in \Int_\phi(R_{\theta_0}): \text{ for all }k\geq 0, \psi^k(y) \in R_{\theta_k}\right\} \\
        W_{loc}^u(\theta) & := \left\{y \in \Int_\phi(R_{\theta_0}): \text{ for all }k\geq 0, \psi^{-k}(y) \in R_{\theta_{-k}}\right\}.
    \end{align*}
\end{definition}

There is a  canonically defined map $\pi: \Theta(\cR) \to \cR$ given by $$\pi(\theta) = W_{loc}^s(\theta) \cap W_{loc}^u(\theta) = \bigcap_{j=-\infty}^\infty \psi^{-j}(R_{\theta_j}).$$

Since the flow is expansive (Proposition \ref{theorem:anosovflowexpansive}), this intersection is a single point on $R_{\theta_0}$ and the map $\pi$ is well-defined. On the other hand, for any point $y \in \cR$, we have a canonical sequence $\theta^y \in \Theta(\cR)$ given by
$$\theta^y = (\ldots, \theta^y_{-1} \mid \theta^y_{0}, \theta^y_{1}, \ldots ) \: \text{ where } \psi^k(y) \in R_{\theta^y_k}.$$
Therefore, $\pi(\theta^y) = y$. 

\subsection{Sections for Flow on $GX$}
In this section, we utilize attributes of the metric Anosov flow to define the tools and verify the properties needed to construct a ``discretized'' subaction. Let $(X, d_X)$ be a locally CAT(-1) space with geodesic flow $\{g_t\}$ which is metric Anosov (Theorem \ref{theorem:GeodesicFlowMetricAnosov}). 

Let $\Sigma$ be the disjoint union of the rectangle sections $\{\Sigma_i\}_{i\in I}$ in the Markov proper family given by Proposition \ref{prop:MarkovProperFamiliesExist}. Let $\psi: \Sigma \to \Sigma$ be the Poincare return map and let $\tau: \Sigma \to (0, \infty)$ be the return time map with respect to the flow $\{g_t\}$. 

While constructing collections of sections at arbitrarily small scales is straightforward (see Proposition 3.4), the key aspect of the next lemma is the ability to construct collections whose sections have arbitrarily small diameters and such that the associated return time is bounded below by a fixed constant, $\tau_*$. 

\begin{lemma}
[cf.  Lemma 6, \cite{LT}]\label{lemma:SectionsCanBeAsSmallAsDesired}
There exist constants $\alpha^*, \tau_*, \tau^*$ and a collection of Poincare sections $\Sigmat=\{\Sigmat_i\}_{i\in I}$ with $\diam\: \Sigmat_i \leq \alpha^*$ such that for any $\alpha < \alpha^*$,  we can construct a collection of disjoint Poincare sub-sections $\Sigma=\{\Sigma_i\}_{i\in I}$ satisfying the following properties:
\begin{enumerate}
    \item for all $i \in I$,  $\Sigma_i \subset \Sigmat_{i}$, and $\diam\: \Sigma_i \leq \alpha$; 
    \item $\bigcup_{i\in I} g_{(0, \tau^*)}\Int_{g} \Sigma_i = GX$;
    \item the return time $\tau$ associated to the sub-section $\Sigma$ is $\geq \tau_*$.
\end{enumerate}
\end{lemma}

\begin{proof}
    We first start with any collection of Poincare sections $\Sigmah = \{\Sigmah_i\}_{i \in \hat I}$, given  by Proposition \ref{prop:MarkovProperFamiliesExist} and suppose the associated return time is bounded below and above by $t_*$ and $t^*$, respectively. Note that $t_*>0$ since the sections in $\Sigmah$ are disjoint and closed. By definition, the $(0, t^*)$-flow boxes $g_{(0, t^*)}\Int_{g} \Sigmah_i$ of $\Sigmah_i$ cover $GX$; that is, $GX = \bigcup_{i\in I} g_{(0, t^*)}\Int_{g} \Sigmah_i$. Choose $\alpha^*$ small enough that $2\alpha^* < t_*$ and the $(0, t^* - 2\alpha^*)$-flow boxes of $\Sigmah'_i$ still cover $GX$, where
    $$\Sigmah'_i = \left\{\gamma \in \Sigmah_i: d(\gamma, GX \backslash \Sigmah_i) > 2\alpha^*\right\}.$$
    Let $C$ be the covering constant (depending only on $X$) given by Lemma \ref{lem:covering}. Note that each $\Sigmah_i'$ can be identified with a subset of $\partial^{(2)}_\infty\tilde X$ because the parametrization of these geodesics is fixed by virtue of their belonging to $\Sigmah'_i$. Hence we can apply Lemma \ref{lem:covering} to it. Given any $\alpha< \alpha^*$, we can cover each $\Sigmah'_i$ with sets $\{B_{ijk}\}_{j \in J, \: 1\leq k \leq C}$ in such a way that for fixed $i, k,$ the sets $\{B_{ijk}\}_{j\in J}$ are pairwise disjoint, and (using Lemma \ref{lem:geo flow is Lipschitz}) such that for any $t\in[0,\alpha^*]$, $\diam\: g_tB_{ijk}<\alpha$. We note here that in the proof of Lemma \ref{lem:covering} the sets $B_{ijk}$ are $\epsilon$-balls with respect to the $\ell^\infty$ product of visual metrics. It is not hard to see that by taking $\epsilon$ small enough, we can ensure that the resulting sets have $d_{GX}$-diameter as small as we need.

    We then `stack' $C$ copies of $\Sigmah_i$ along the flow, letting
    $$\Sigmat_{ijk} = g_{\frac{k}{C}\alpha^*} \Sigmah_{i} \: \text{ and } \: \Sigma_{ijk} = g_{\frac{k}{C}\alpha^*} B_{ijk}.$$
    Note that by the choice of $\alpha^*$, $\Sigma_{ijk} \subset \Sigmat_{ijk}$.
    Let $\tau$ denote the return time associated to the sections $\Sigma = \{ \Sigma_{ijk}\}_{i,j,k}$.

    By design, for fixed $i$, $\{\Sigma_{ijk}\}_{k, j}$ are pairwise disjoint, but additionally, since the maximum height of the stack is $\alpha^* < \frac{t_*}{2}$, no two stacks intersect; therefore, $\{\Sigma_{ijk} \}_{i,j,k}$ are all pairwise disjoint. Also, by construction, $\tau$ is at least $\frac{\alpha^*}{C} =: \tau_*$. Finally, we note that the longest $\tau$ can be  is the sum of the maximum return time for $\Sigmah'$ and the maximum stack height; that is $\tau < (t^* -2\alpha^*) + \alpha^* = t^* -\alpha^* =: \tau^*$, and so the $(0, \tau^*)$-flow boxes for $\Sigma$ cover $GX$:
    $$\bigcup_{i,j,k} g_{(0, \tau^*)} \Int_{g} \Sigma_{ijk} \supset \bigcup_{i} g_{(0, t^*-2\alpha^*)} \Int_{g} \Sigmah'_i = GX.$$
    Importantly, we note that $\tau_*$ and $\tau^*$ are independent of the size $\alpha$ of the sections $\Sigma$. 
\end{proof}

\begin{figure}[h!]
\centerline{\includegraphics[width=1.3\textwidth]{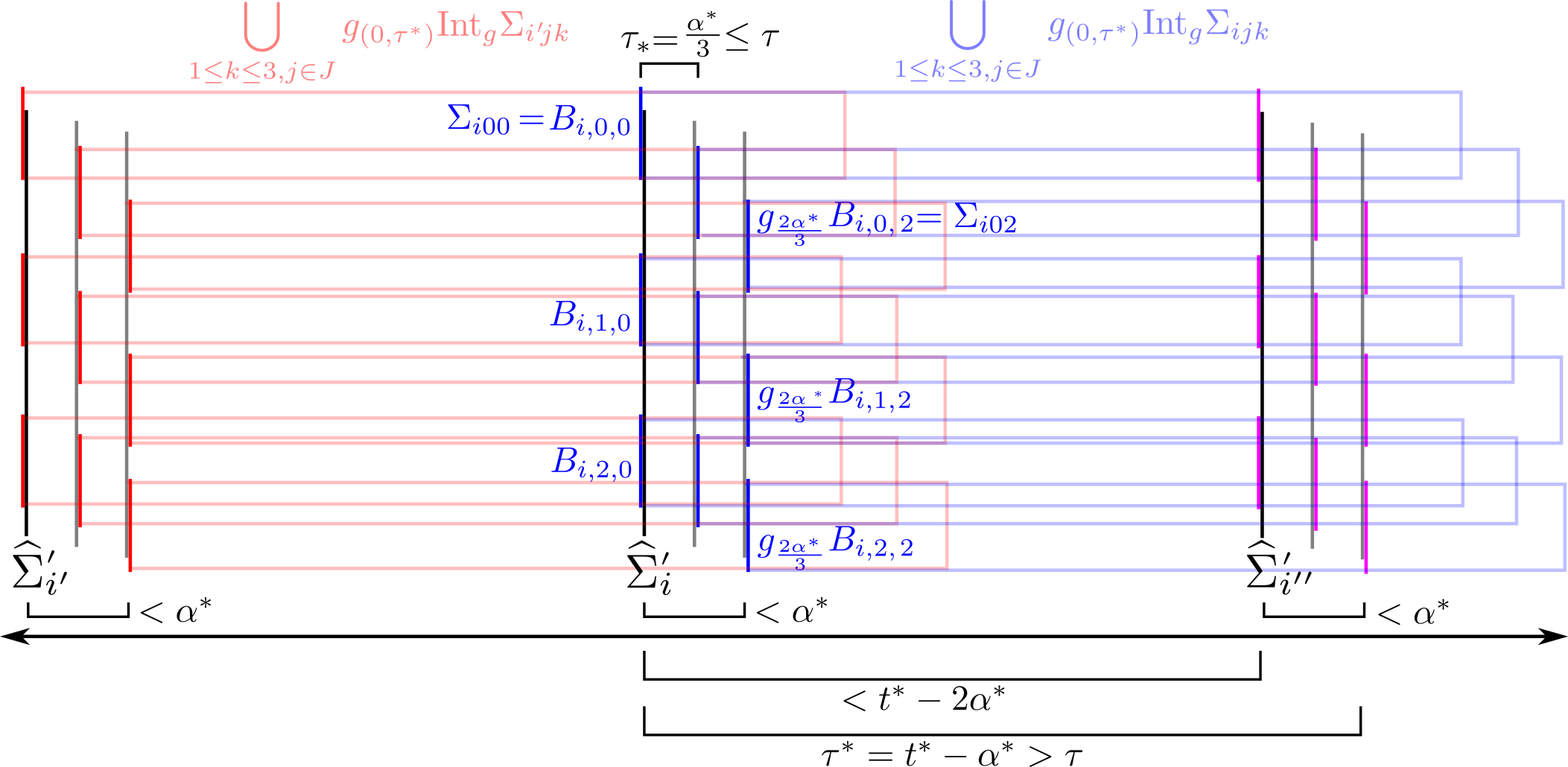}}
	\caption{A visual aid for the proof of Lemma \ref{lemma:SectionsCanBeAsSmallAsDesired}. Note that the stack associated to $\Sigma_{ijk}$ is covered by the flow boxes for $\Sigma_{i'jk}$ (shown in light red). Here, $C = 3$.}
	\label{fig:covering}
\end{figure}

From now on, we fix the sections $\Sigmat$ with fixed associated constants $\alpha^*, \tau_*, \tau^*$ given by Lemma \ref{lemma:SectionsCanBeAsSmallAsDesired}; in particular, $\diam \:\Sigmat_i \leq \alpha^*$ and its associated return time $\taut$ is bounded below and above by the fixed constants $\tau_*$ and $\tau^*$, respectively. By the lemma, we will be able to choose subsections $\Sigma$ such that $\diam\:\Sigma_i= \alpha$ is arbitrarily small, while its associated return time $\tau$ preserves the same bounds.

\begin{definition}
    Let $i, j \in I$. We say $i \to j$ is a \textit{simple transition} if there exists $\gamma \in \Sigma_i$ such that $\psi(\gamma) \in \Sigma_j$. Let
    $$\dom(\psit_{ij})=\dom(\taut_{ij}) := \{\gamma \in \tilde{\Sigma}_{i}: g_t(\gamma) \in \tilde{\Sigma}_{j} \text{ for some } t \in  (0, \tau^*) \}.$$
    We define the extended first return time $\taut_{ij}$ and first return  map $\tilde{\psi}_{ij}$ for the simple transition $i \to j$ as follows: for $\gamma \in \dom(\psit_{ij})=\dom(\taut_{ij})$, 
    $$\taut_{ij}(\gamma) = \inf\{t \in (0, \tau^*): g_t\gamma \in \tilde{\Sigma}_{j}\}$$
    $$\psit_{ij}(\gamma) = g_{\taut_{ij}(\gamma)}\gamma.$$
\end{definition}
In particular, using Lemma \ref{lemma:SectionsCanBeAsSmallAsDesired}, we choose $\alpha$ to be small enough that $\Sigma_i \subseteq \dom(\psit_{ij})$ and $\Sigma_j \subseteq \range(\psit_{ij})$ for any simple transition $i \to j$.

Note that for any canonical sequence $\theta \in \Theta(\Sigma)$ (see Definition \ref{def:CanonicalCodingSpace}), we have $\theta_{k} \to \theta_{k+1}$ is a simple transition for each $k \in \mathbb{Z}$. Indeed, 
$\psi\big(\psi^{k}(\pi(\theta)) \big) = \psi^{k+1}(\pi(\theta)) \in \Sigma_{\theta_{k+1}}$. We now extend the canonical coding space $\Theta(\Sigma)$ to allow for double-sided sequences of simple transitions. 

\begin{definition}
    A \textit{pseudo-orbit} is a double-sided sequence of simple transitions. Let $\Omega$ be the extended coding space of pseudo-orbits
    $$\Omega = \{\omega = (\cdots, \omega_{-1} \mid \omega_{0}, \omega_1, \cdots): \: \omega_k \to \omega_{k+1} \text{ is a simple transition for all } k\}.$$
\end{definition}
It is immediate from its definition that $\Omega$ is also a sub-shift of finite type, and we denote by $\sigma: \Omega \to \Omega$ the associated left shift map. Note that $\Theta \subset \Omega$. For each $\omega \in \Omega$, define 
$$\psit_\omega = \psit_{\omega_0 \omega_1}: \Sigma_{\omega_0} \to \Sigmat_{\omega_1}$$
$$\taut_\omega = \taut_{\omega_0 \omega_1}$$
and 
$$\psit_\omega^k = \psit_{\sigma^{k-1}(\omega)}\circ \cdots \circ \psit_{\omega}.$$

As before, we define the local stable and unstable sets of a pseudo-orbit $\omega \in \Omega$:
\begin{definition}\label{def:stableunstableforomega}
     For any $\omega \in \Omega$, we define the local stable and unstable sets of $\omega$ to be:
    \begin{align*}
        W_{loc}^s(\omega) & := \left\{\gamma \in \Int_\phi(\Sigma_{\omega_0}): \text{ for all }k\geq 0, \psit_{\omega}^k(\gamma) \in \Sigma_{\omega_k}\right\} \\
        W_{loc}^u(\omega) & := \left\{\gamma \in \Int_\phi(\Sigma_{\omega_0}): \text{ for all }k\geq 0, \psit_{\omega}^{-k}(\gamma) \in \Sigma_{\omega_{-k}}\right\}.
    \end{align*}
\end{definition}

For any $\omega, \omega' \in \Omega$ with $\omega_0 = \omega_0'$, $W_{loc}^s(\omega)$ intersects $W_{loc}^u(\omega')$ at a unique point. This is due to the fact that $\alpha$ can be made arbitrarily small and Bowen's shadowing lemma: consider the sequence formed by concatenating the past of $\omega'$ with the future of $\omega$: $\zeta = (\cdots, \omega'_{-2}, \omega'_{-1} \mid \omega_0, \omega_1, \cdots)$. By definition, for each $j$,  $\zeta_j \to \zeta_{j+1},$ so there exists $\gamma_j \in \Sigma_{\zeta_j}$ such that $\psi_\zeta(\gamma_j) \in \Sigma_{\zeta_{j+1}}$. Therefore, $(\gamma_j)_{j \in \Z}$ is a pseudo-orbit and by Bowen's shadowing lemma, there exists a unique element $\gamma$ whose true orbit $(\psi_\zeta(\gamma))_{j \in \Z}$ shadows $(\gamma_j)_{j \in \Z}$. Therefore, $\{\gamma\} = W_{loc}^s(\omega) \cap W_{loc}^u(\omega')$.  

\begin{definition}
For any $\omega, \omega' \in \Omega$ with $\omega_0 = \omega_0'$, we denote the unique point in the intersection $W_{loc}^s(\omega) \cap W_{loc}^u(\omega')$ by $[\omega, \omega']$. We denote $[\omega, \omega]$ by $\pi(\omega)$.      
\end{definition}

We now demonstrate that $\psit_\omega$ satisfies a property that is related to and resembles the Markov Property (Definition \ref{def:Markov proper family}).

\begin{lemma}[Markov Property for $\psit_\omega$]\label{lem:markov} $\psit_\omega$ stabilizes the local stable set:
    $$\psit_\omega(W^s_{loc}(\omega)) \subset W^s_{loc}(\sigma(\omega))$$
and $\psit^{-1}_\omega$ stabilizes the local unstable set:
    $$\psit_\omega^{-1}(W^u_{loc}(\omega)) \subset W^u_{loc}(\sigma^{-1}(\omega)).$$
Consequently, $\psit_\omega(\pi(\omega)) = \pi(\sigma(\omega))$. 
\end{lemma}

\begin{proof}
    Let $\gamma \in W^s_{loc}(\omega)$. Then for any $k \geq 0$, $$\psit_{\sigma(\omega)}^k\big(\psit_\omega(\gamma)\big) = \psit_\omega^{k+1}(\gamma) \in \Sigma_{\omega_{k+1}} = \Sigma_{(\sigma(\omega))_{k}}$$
    so $\psit_\omega(\gamma) \in W^s_{loc}(\sigma(\omega))$. Similarly, let $\gamma \in W^u_{loc}(\omega)$. Then for any $k \geq 0$, $$\psit_{\sigma^{-1}(\omega)}^{-k}\big(\psit^{-1}_\omega(\gamma)\big) = \psit_\omega^{-(k+1)}(\gamma) \in \Sigma_{\omega_{-(k+1)}} = \Sigma_{(\sigma^{-1}(\omega))_{-k}}$$
    so $\psit^{-1}_\omega(\gamma) \in W^u_{loc}(\sigma^{-1}(\omega))$. 
\end{proof}

\subsection{Exponential contraction}
\label{sec:lemma 8 replacement}

Exponential contraction along $W^s_{loc}(\theta)$ under the return map $\psi$ for any $\theta\in \Theta$ plays an essential role in later proofs. We will prove our contraction result (an analogue for our situation of Lemma 8(iii) in \cite{LT}) for the simple first return map $\psi$ associated to the collection $\Sigma$ if sections constructed above. We note that the result holds for the `return map with instructions' $\psit_\omega$ since $\psit_\omega$ is at any point $\psi^k$ for an appropriate $k$.

\begin{lemma}\label{lem:exp contraction}
    Let $\gamma, \gamma'\in W^s_{loc}(\omega)$ for some $\omega\in \Omega$ and let $\psi$ be the first return map for the collection of sections $\Sigma$. Then there exists a constant $C>0$ (independent of $\gamma,\gamma',\omega$) such that for all $k>0$,
    \[d_{GX}(\psi^k \gamma, \psi^k \gamma') < Ce^{-k\tau_*} d_{GX}(\gamma,\gamma').\]
    Similarly, if $\gamma,\gamma'\in W^u_{loc}(\omega)$, for all $k>0$
    \[d_{GX}(\psi^{-k} \gamma, \psi^{-k} \gamma') < Ce^{-k\tau_*} d_{GX}(\gamma,\gamma').\]    
\end{lemma}

\begin{remark}
    For comparison with \cite{LT} in this and subsequent lemmas, recall that $\lambda=1$ for geodesic flow on a CAT(-1) space.
\end{remark}

\begin{proof}
We prove the result only for the local stable sets, the proof for the unstable sets being essentially the same. As usual, we lift to the CAT(-1) universal cover to make our geometric arguments.

Suppose that $\gamma, \gamma'\in W^s_{loc}(\omega) \subset \Sigma_i$ and that $\psi^k\gamma, \psi^k \gamma' \in \Sigma_j$. 
Since $\gamma, \gamma'\in W^s_{loc}(\omega)$, $\gamma'$ is on the weak-stable leaf through $\gamma$ and hence for some $r_1$, $g_{r_1}\gamma'\in W^{ss}(\gamma)$. Similarly, as $\psi^k \gamma, \psi^k \gamma'\in W^s_{loc}(\sigma^k\omega)$, for some $r_2$, $g_{r_2} \psi^k \gamma' \in W^{ss}_{loc}(\psi^k \gamma)$. Let $\tau, \tau'$ be such that $g_\tau \gamma = \psi^k \gamma$ and $g_{\tau'} \gamma' = \psi^k \gamma'$. (See Figure \ref{fig:exp contraction}.) Note that by construction of $\Sigma$, $\tau, \tau' \geq k \tau_*$.

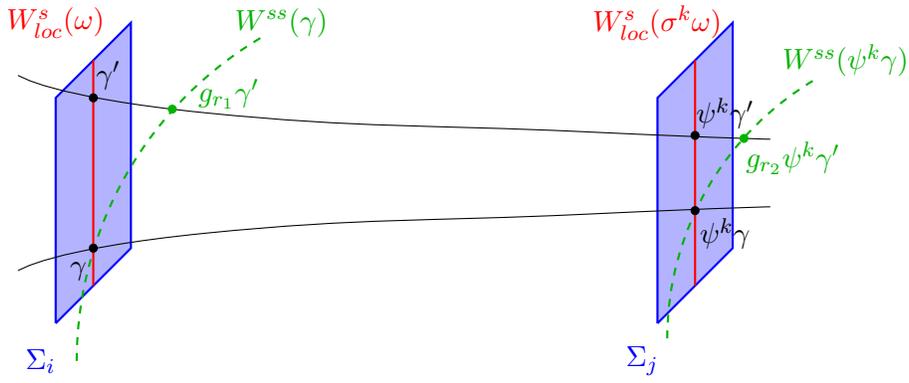
\begin{figure}[h]
\centering
\begin{tikzpicture}[scale=1]

\draw[thick, blue] (-4,-2) -- (-4,1) -- (-3,2) -- (-3,-1) -- (-4,-2);
\fill[blue, opacity=.3] (-4,-2) -- (-4,1) -- (-3,2) -- (-3,-1) -- (-4,-2);
\draw[thick,blue] (4,-2) -- (4,1) -- (5,2) -- (5,-1) -- (4,-2);
\fill[blue, opacity=.3] (4,-2) -- (4,1) -- (5,2) -- (5,-1) -- (4,-2);

\draw[thick, red] (-3.5,-1.5) -- (-3.5,1.5) ;
\draw[thick, red] (4.5,-1.5) -- (4.5,1.5) ;

\draw [black] plot [smooth, tension=1] coordinates { (-4.5,-1.3) (-2,-.8) (3,-.55) (5.5,-.45)};
\draw [black] plot [smooth, tension=1] coordinates { (-4.5,1.3) (-2,.8) (3,.55) (5.5,.45)};

\draw[thick, green!70!black, dashed] (-3.72,-2.5) arc (180:120:5);
\draw[thick, green!70!black, dashed] (4.13,-2.2) arc (180:120:4);

\fill (-3.5,-1) circle(.06);
\fill (-3.5,1) circle(.06);
\fill (4.5,-.5) circle(.06);
\fill (4.5,.5) circle(.06);

\fill[green!70!black] (-2.45,.85) circle(.06);
\fill[green!70!black] (5.15,.46) circle(.06);

\node at (-3.3, 1.3) {$\gamma'$} ;
\node at (-3.7, -1.3) {$\gamma$} ;
\node at (4.9, .8) {$\psi^k \gamma'$} ;
\node at (4.9, -.8) {$\psi^k \gamma$} ;
\node[blue] at (-4.2, -2.5) {$\Sigma_i$} ;
\node[blue] at (3.8, -2.5) {$\Sigma_j$} ;
\node[red] at (-4, 2) {$W^s_{loc}(\omega)$} ;
\node[red] at (4, 2) {$W^s_{loc}(\sigma^k\omega)$} ;
\node[green!70!black] at (-1.7,1.05) {$g_{r_1}\gamma'$} ;
\node[green!70!black] at (-1,2) {$W^{ss}(\gamma)$} ;
\node[green!70!black] at (5.8,.2) {$g_{r_2}\psi^k \gamma'$} ;
\node[green!70!black] at (6.5,1.5) {$W^{ss}(\psi^k \gamma)$} ;

\end{tikzpicture}
\caption{The geometric setup for Lemma \ref{lem:exp contraction}.}\label{fig:exp contraction}
\end{figure}

Since $g_\tau$ maps $W^{ss}(\gamma)$ to $W^{ss}(\psi^k \gamma)$ we see that $g_{r_2} \psi \gamma' = g_\tau g_{r_1} \gamma'$ and so $\tau'=\tau+(r_1-r_2).$

By the triangle inequality, 
\[d_{GX}(\gamma, g_{r_1}\gamma') \leq d_{GX}(\gamma,\gamma') + |r_1|,\]
\[d_{GX}(\psi^k\gamma, \psi^k\gamma') \leq d_{GX}(\psi^k \gamma, g_{r_2}\psi^k \gamma') +|r_2|.\]
Since $g_t$ is a metric Anosov flow with $\lambda=1$, for some uniform $C_1>0$
\[d_{GX}(\psi^k \gamma, g_{r_2}\psi^k \gamma') \leq C_1 e^{-k \tau_*} d_{GX}(\gamma, g_{r_1}\gamma').\]
Combining these inequalities gives
\begin{align}
    d_{GX}(\psi^k \gamma, \psi^k \gamma') &\leq d_{GX}(\psi^k \gamma, g_{r_2}\psi^k \gamma') + |r_2| \nonumber \\
        & \leq C_1 e^{-k\tau_*} d_{GX}(\gamma, g_{r_1}\gamma')+|r_2| \nonumber \\
        & \leq C_1 e^{-k\tau_*} [d_{GX}(\gamma,\gamma')+|r_1|]+|r_2|. \nonumber
\end{align}

We claim the following sublemma:

\begin{sublemma}
$|r_1|\leq 2 d_{GX}(\gamma,\gamma')$ and $|r_2|\leq 2 d_{GX}(\psi^k \gamma, g_{r_2}\psi^k \gamma')$.
\end{sublemma}

Assuming the sublemma,
\begin{align}
    d_{GX}(\psi^k \gamma, \psi^k \gamma') &\leq C_1 e^{-k \tau_*}[d_{GX}(\gamma,\gamma')+ 2 d_{GX}(\gamma,\gamma')]+2d_{GX}(\psi^k \gamma, g_{r_2}\psi^k \gamma') \nonumber \\
        & \leq 3C_1 e^{-k\tau_*} d_{GX}(\gamma,\gamma') + 2C_1 e^{-k\lambda\tau_*} d_{GX}(\gamma, g_{r_1}\gamma') \nonumber \\
        & \leq 3C_1 e^{-k\tau_*} d_{GX}(\gamma,\gamma') + 2C_1 e^{-k\tau_*}[d_{GX}(\gamma,\gamma')+r_1] \nonumber \\
        & \leq 3C_1 e^{-k\tau_*} d_{GX}(\gamma,\gamma') + 2C_1 e^{-k\tau_*}[d_{GX}(\gamma,\gamma') + 2d_{GX}(\gamma,\gamma')] \nonumber \\
        & = 9C_1 e^{-k\tau_*} d_{GX}(\gamma,\gamma') \nonumber
\end{align}
as desired.

\begin{proof}[Proof of Sublemma]
We prove the sublemma for $r_1$; the proof for $r_2$ is the same.

Note that since the strong stable set $W^{ss}(\gamma)$ is determined by values of $B_{\gamma(0)}(-,\gamma(+\infty))$, $|r_1| = B_{\gamma(0)}(\gamma'(0),\gamma(+\infty))$. Busemann functions are 1-Lipschitz, so $B_{\gamma(0)}(\gamma'(0),\gamma(+\infty)) \leq d_X(\gamma(0), \gamma'(0))$. Therefore
\[|r_1|\leq d_X(\gamma(0), \gamma'(0)) \leq 2d_{GX}(\gamma,\gamma')\]
as desired.
\end{proof}
\end{proof}

We also have the following analogue of Lemma 8(iv) from \cite{LT} which gives very rough upper and lower bounds on the distance between two points under iterates of our return map $\tilde \psi$.

\begin{lemma}\label{lem:bounds on d(psix, psiy)}
There exist constants $K^*$ and $\Lambda_*^s < 0 < \Lambda_*^u$ such that the following holds.
Suppose that $x,y\in \Sigma_i$ and that $i=i_0 \to \cdots \to i_n$ is a chain of simple transitions so that $\tilde \psi^n:=\tilde \psi_{i_{n-1}i_n}\circ \cdots \circ \tilde \psi_{i_0i_1}$ exists at $x$ and $y$. Then
\[(K^*)^{-1}\exp(n\Lambda_*^s)d_{GX}(x,y) \leq d_{GX}(\tilde \psi^n x,\tilde \psi^n y) \leq K^* \exp(n\Lambda_*^u)d_{GX}(x,y).\]
\end{lemma}

\begin{proof}
Say $\tilde \psi^n x = g_tx$ and $\tilde \psi^n y = g_{t'}y$. Note that $t\leq n\tau^*$ and that $|t-t'|$ is the difference in first return times for $x$ and $y$ and under flow from $\Sigma_{i}$ to $\Sigma_{i_n}$.

By triangle inequality, 
\[d_{GX}(\tilde \psi^n x, \tilde \psi^n y) \leq d_{GX}(g_tx, g_ty)+ |t-t'|.\]
By Lemma \ref{lem:geo flow is Lipschitz} and the Lipschitz behavior of return times given by Proposition \ref{prop:MarkovProperFamiliesExist}, 
\begin{align}
    d_{GX}(\tilde \psi^n x, \tilde \psi^n y) &\leq e^{2t}d_{GX}(x, y)+ Kd_{GX}(x, y) \nonumber \\
        & \leq (K+1) \exp((2\tau^*)n)d_{GX}(x,y).\nonumber
\end{align}

To prove the other bound we can simply reverse the direction of the flow: $x = g_{-t}\tilde \psi x$ and $y = g_{-t'}\tilde \psi y$. Applying the same argument we get
\[d_{GX}(x,y)\leq (K+1) \exp((2\tau^*)n)d_{GX}(\tilde \psi^nx, \tilde \psi^ny)\]
completing the proof.
\end{proof}

We now prove two lemmas that will be useful in proving that the discretized subaction is H\"older in Section \ref{sec:discretized}.

\begin{lemma}\label{lem:Lemma 23 replacement}
    There exists $\delta > 0$ depending only on $\Sigma$ such that for any $i \in I$ and any $\gamma, \gamma' \in \Sigma_i$, if $d_{GX}(\gamma, \gamma') < \delta$, then there exists a simple transition $i \to j$ and $m, n \geq 1$ such that 
    $$\psi^m(\gamma) \in \Sigma_j,\: \psi^n(\gamma') \in \Sigma_j $$ with first return times from $\Sigma_i$ to $\Sigma_j$
    $$\tau_m(\gamma) \leq \tau^*, \: \tau_n(\gamma') \leq \tau^*.$$
\end{lemma}

\begin{proof}
    Let $U_i = g_{(0, \tau^*)} \Sigma_i$. For any simple transition $i \to j$, there is a map $t_j: \Delta_{ij} \to (0, \tau^*)$ such that 
    $$\Sigma_j \cap U_i = \{g_{t_j(\eta)}(\eta): \eta \in \Delta_{ij} \}.$$
    Since $\tau^*$ is the maximum return time, the sets $\{\Delta_{ij}\}_j$ cover $\Sigma_i$, which is compact. Therefore, we can find a Lebesgue number $\delta > 0$ for the cover $\{\Delta_{ij}\}_j$. If $\gamma'$ is contained in the $\delta$-ball around $\gamma$, then there exists $j\in I$ such that the ball is contained in $\Delta_{ij}$. In particular, $\gamma, \gamma' \in \Delta_{ij}$ with $t_j(\gamma) =: \tau_m(\gamma) \leq \tau^*$ and $t_j(\gamma') =: \tau_n(\gamma') \leq \tau^*$, as desired.
\end{proof}

In the following, $\Lambda_*^u$ and $K^*$ are constants given by Lemma \ref{lem:bounds on d(psix, psiy)}.

\begin{lemma}\label{lem:lemma24replacement}
    Suppose $\gamma, \gamma' \in \Sigma_i$. For any $N \geq 1$, if $d_{GX}(\gamma, \gamma') < \frac{\delta}{K^*} e^{-N\Lambda_*^u}$, then there exist $\omega, \omega' \in \Omega$ such that $\pi(\omega) = \gamma$, $\pi(\omega') = \gamma'$ and their first $N+1$ symbols coincide:
    $\omega_0= \omega'_0, \ldots, \omega_N = \omega'_N.$
\end{lemma}

\begin{proof}
Let $\gamma_0= \gamma,$ and $\gamma'_0= \gamma' \in \Sigma_i$. By Lemma \ref{lem:Lemma 23 replacement}, since $d_{GX}(\gamma_0, \gamma'_0) < \delta$, there exists $m_1, n_1$ and $i_1 \in I$ such that 
$$\gamma_1 = \psi^{m_1}(\gamma_0), \: \gamma'_1 = \psi^{n_1}(\gamma'_0) \in \Sigma_{i_1}.$$
By Lemma \ref{lem:bounds on d(psix, psiy)}, 
\begin{align*}
    d_{GX}(\gamma_1, \gamma'_1) \leq K^*  e^{\Lambda_*^u} d_{GX}(\gamma_0, \gamma'_0) < K^* e^{\Lambda_*^u} \frac{\delta}{K^*} e^{-N\Lambda_*^u} \: = \delta e^{-(N-1)\Lambda_*^u}\:<\delta
\end{align*}
so we can apply Lemma \ref{lem:Lemma 23 replacement} again. We repeat this construction $N$ times to find two sequences $(m_1, m_2, \ldots, m_N)$ and $(n_1, \ldots, n_N)$ such that for each $0 \leq l \leq N$
$$\gamma_l = \psi^{m_1 + \ldots + m_l}(\gamma), \: \text{ and } \gamma'_l = \psi^{n_1 + \ldots + n_l}(\gamma')$$
belong to the same section $\Sigma_{i_l}$ with $d_{GX}(\gamma_l, \gamma'_l) < \delta$.

Let $\theta^\gamma, \theta^{\gamma'}\in \Theta(\Sigma)$ be the canonical sequences associated to $\gamma$ and $\gamma'$, respectively, and let 
$$m = m_1 + \cdots + m_N, \: \: n= n_1 + \cdots + n_N.$$
Define
\begin{align*}
    \omega & = (\ldots, \theta_{-2}^{\gamma}, \theta_{-1}^{\gamma} \mid i_0, i_1, \ldots, i_N, \theta_{m+1}^{\gamma}, \theta_{m+2}^{\gamma}, \ldots) \\
    \omega' & = (\ldots, \theta_{-2}^{\gamma'}, \theta_{-1}^{\gamma'} \mid i_0, i_1, \ldots, i_N, \theta_{n+1}^{\gamma'}, \theta_{n+2}^{\gamma'}, \ldots)
\end{align*}
By construction, $\psit_{\omega}^k(\gamma) \in \Sigma_{\omega_k}$ and $\psit_{\omega'}^k(\gamma') \in \Sigma_{\omega'_k}$ for all $k \in \Z$, so $\gamma = \pi(\omega)$ and $\gamma' = \pi(\omega')$, and they coincide in the first $N+1$ times, as desired.
\end{proof}

\section{A sub-action for the discretized system}\label{sec:discretized}

Let $A: GX \to \R$ be  H\"older and define its \emph{minimal average} by
$$m(A) = \inf\left\{ \int A \; d\mu \mid \mu \in \mathcal{M}_1(GX, g_t)\right \}$$
where $\mathcal{M}_1(GX, g_t)$ denotes the set of all $g_t$-invariant probability measures on $GX$. 
We define the discretized observable $\mathcal{A}: \Sigma \to \R$ by
\begin{equation}
    \mathcal{A}(\gamma) = \int_0^{\tau(\gamma)} (A - m(A)) \circ g_t (\gamma) \; dt.
\end{equation}
The goal of this section is to construct a discretized sub-action $\mathcal{V}$ for $\mathcal{A}$:

\begin{proposition}\label{prop:DiscretizedSubactionV}
    There exist sections $\Sigma$ and a globally H\"older function $\mathcal{V}:\Sigma \to \R$ such that
    $$\mathcal{A}(\gamma) \geq \mathcal{V} \circ \psi(\gamma) - \mathcal{V}(\gamma) \: \text{for all } \gamma \in \Sigma.$$
\end{proposition}

As in \cite{LT}, we first extend $\mathcal{A}: \Sigma \to \R$ to the space $\Omega \times \Sigma$. To that end, we introduce the following notation. Define
\begin{align*}
    \psit: \Omega \times \Sigma & \to \Omega \times \Sigma \\
    (\omega, \gamma) & \mapsto (\sigma(\omega), \psit_\omega(\gamma)),
\end{align*}
to be the first return map ``with instructions'' and let $\taut$ be the associated first return time $\taut(\omega, \gamma) := \taut_{\omega}(\gamma)$. We also extend the discretized $\mathcal{A}$ to $\Omega \times \Sigma$ by:
$$\tilde{\mathcal{A}}(\omega, \gamma) = \int_0^{\taut(\omega, \gamma)} (A - m(A)) \circ g_t(\gamma) \; dt.$$
Recall that for $\gamma \in \Sigma$, its canonical pseudo-orbit $\theta^\gamma \in \Theta(\Sigma)$ given by
$$\theta^\gamma = (\ldots, \theta_{-1} \mid \theta_{0}, \theta_{1}, \ldots ) \: \text{ where } \psi^k(\gamma) \in \Sigma_{\theta_k}$$
satisfies $\pi(\theta^\gamma)=  \gamma$. Define
\begin{align*}
    \tilde{\theta}: \Sigma  & \to \Omega \times \Sigma \\
   \gamma & \mapsto (\theta^\gamma, \gamma).
\end{align*}
Then we note that
\begin{equation*}
    \tilde{\theta}\circ \psi = \psit \circ \tilde{\theta}.
\end{equation*}
Since $\taut(\theta^\gamma, \gamma) = \tau(\gamma)$, we note that $\tilde{\mathcal{A}}$ is an extension of $\mathcal{A}$:
$$\tilde{\mathcal{A}} \circ \tilde{\theta} = \mathcal{A}.$$

\begin{figure}[h!]
	\centering
	\includegraphics[width=0.8\textwidth]{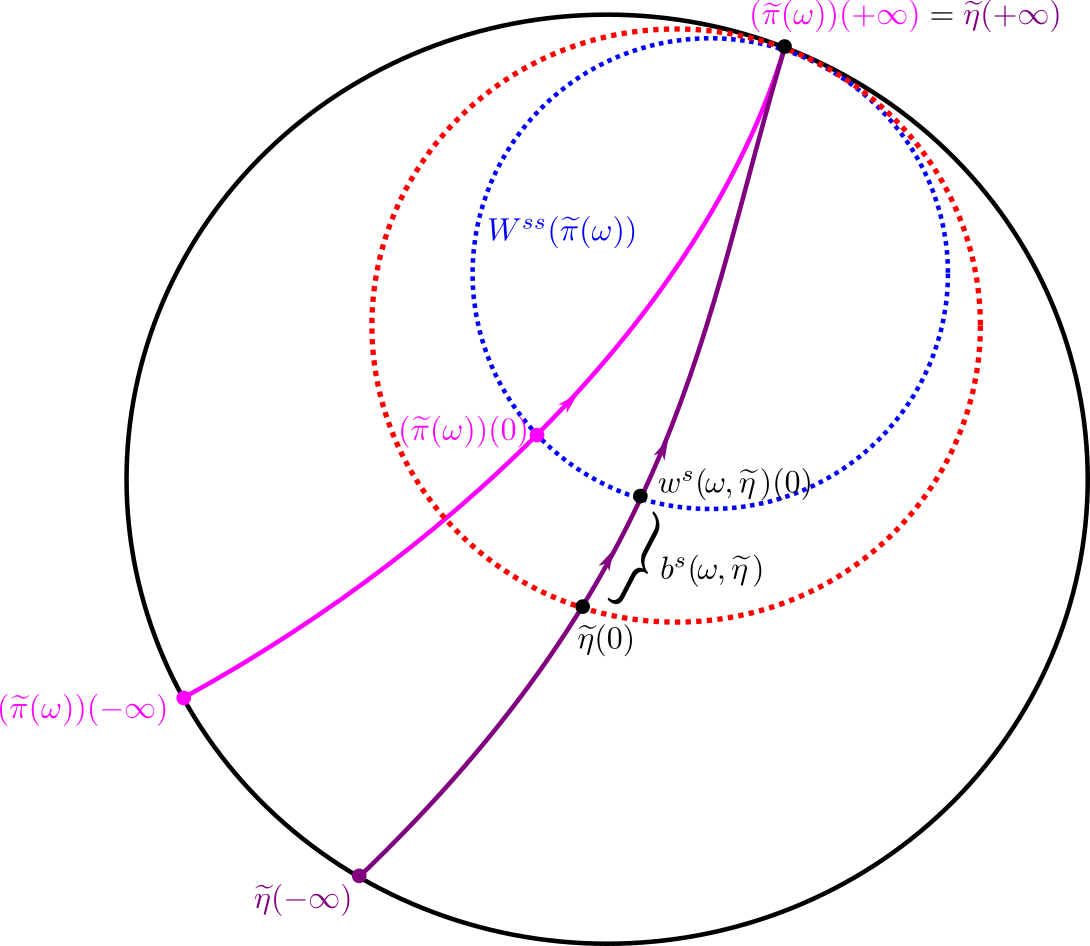}
	\caption{An illustration of some of the definitions from Definition \ref{def:bwdelta}. Note that $\eta \in W^s_{loc}(\omega)$ so $\etat \in W^{cs}(\Tilde{\pi}(\omega))$.}
	\label{fig:busemann}
\end{figure}

\begin{definition}\label{def:bwdelta}
For $\omega \in \Omega$ and $\eta \in W_{loc}^s(\omega) \cap \Sigma_{\omega_0}$, define 
\begin{align*}
    b^s(\omega, \eta) & = B_{\pi(\omega)}(\eta(0))\\
    w^s(\omega, \eta) & = g_{b^s(\omega, \eta)} \eta\\
    \Delta^s(\omega, \eta) & = \sum_{n \geq 0} \left(\mathcal{A} \circ \psit^n(\omega, \eta) - \mathcal{A}\circ \psit^n(\omega,\pi(\omega))\right)
\end{align*}
\end{definition}
\begin{remark}
    In the definition above, we are abusing notation slightly: it should be understood that $B_{\pi(\omega)}(\eta(0))$ is in fact $B_{\Tilde{\pi}(\omega)}(\etat(0))$, where $\Tilde{\pi}(\omega)$ is a fixed lift  of $\pi(\omega)$ to the CAT(-1) universal cover and $\etat \in W^{cs}(\Tilde{\pi}(\omega))$ is the lift of $\eta$ that has the same endpoint at infinity $\Tilde{\pi}(\omega)(+\infty)$ (see Figure \ref{fig:busemann}).
\end{remark}

\begin{lemma}\label{lem:prop12replacement}
    Let $\omega \in \Omega$. 
\begin{enumerate}
    \item[(i)] The map $W^s_{loc}(\omega) \to W^{ss}(\pi(\omega))$ given by $\eta \mapsto w^s(\omega, \eta)$ is a parameterization of $W^{ss}(\pi(\omega))$. 
    \item[(ii)] The stable cocycle $\Delta^s(\omega, \eta)$ admits the equivalent form
    \begin{align*}
        \Delta^s(\omega, \eta) & = \int_0^\infty (A \circ g_t w^s(\omega, \eta) - A \circ g_t \pi(\omega)) \; dt  \\
        & \: \: \: - \int_0^{b^s(\omega, \eta)} (A - m(A))\circ g_t \eta \; dt.
    \end{align*}
\end{enumerate}
In particular, $b^s(\omega, \eta), W^s_{loc}(\omega),$ and $\Delta^s(\omega, \eta)$ depend on $\omega \in \Omega$ only via $\pi(\omega)$. 
\end{lemma}
\begin{proof}
For (i), note that $\eta \in W^{cs}(\pi(\omega))$, so using Lemma \ref{lem:BusemannProperties}(2),
    \begin{align*}
        B_{\pi(\omega)}\big(w^s(\omega, \eta)(0)\big) & =  B_{\pi(\omega)}\big(\eta(b^s(\omega, \eta))\big) 
        \\
        & = B_{\pi(\omega)}(\eta(0)) + (0-b^s(\omega, \eta)) \\
        & = 0.
    \end{align*}
Therefore, $w^s(\omega, \eta) \in W^{ss}(\pi(\omega))$. Conversely, any point on $W^{ss}(\pi(\omega))$ is of the form $g_s(\eta)$ for some $\eta \in \Sigma_{\omega_0}$ and small $s$. Since they are on the same stable leaf, $\omega$ is a valid set of forward instructions that $\eta$ can take. Therefore, $\eta \in W^s_{loc}(\omega)$. 

The proof of (ii) is verbatim the same as \cite{LT} Proposition 12(iii) once the following cocycle equation \eqref{eq:CocycleEquation} has been established in our setting. 

\begin{sublemma}
    For $\eta \in W^s_{loc}(\omega)$, the cocycle equation 
    \begin{equation}\label{eq:CocycleEquation}
        b^s(\omega, \eta) + \tau_n(\pi(\omega)) = b^s(\sigma^n\omega, \psit_\omega^n(\eta)) + \tau_n(\eta)
    \end{equation}
    holds, where $\tau_n(-) = \sum_{k = 0}^{n-1} \tau\circ \psit^k_{\omega}(-)$.
\end{sublemma}
\begin{proof}[Proof of Sublemma]
Using Lemma \ref{lem:BusemannProperties}(1) and (2),
    \begin{align*}
        b^s(\omega, \eta) - b^s(\sigma^n\omega, \psit_\omega^n(\eta)) &  = B_{\pi(\omega)}(\eta(0)) - B_{\pi(\sigma^n(\omega))}\big((\psit^n_{\omega}\eta)(0)\big)\\
        & = B_{\pi(\omega)}(\eta(0)) - B_{g_{\tau_n(\pi(\omega))}\pi(\omega)}\big((\psit^n_{\omega}\eta)(0)\big) \\
        & = B_{\pi(\omega)}(\eta(0))- \left[B_{\pi(\omega)}\big((\psit^n_{\omega}\eta)(0)\big) + \tau_n(\pi(\omega)) \right] \\
        & = \left[B_{\pi(\omega)}(\eta(0)) - B_{\pi(\omega)}\big(\eta(\tau_n(\eta))\big)\right] - \tau_n(\pi(\omega)) \\
        & = \tau_n(\eta) - \tau_n(\pi(\omega)).
    \end{align*}
\end{proof}

Finally, if $\pi(\omega) = \pi(\omega')$, then by definition, for any $\eta \in W^s_{loc}(\omega) \cap W^s_{loc}(\omega')$, $b^s(\omega, \eta) = b^s(\omega', \eta)$. Then, (i) gives us that $W^s_{loc}(\omega)= W^s_{loc}(\omega')$, and (ii) gives us that $\Delta^s(\omega, \eta) = \Delta^s(\omega', \eta)$.
\end{proof}

We now define a discretized sub-action $\Vt$ on the space $\Omega$. Let $S_n \At = \sum_{k=0}^{n-1} \At \circ \psit^k$ be the Birkhoff sum of $\At$. 

\begin{definition}
   For any $\omega \in \Omega,$ define
   \begin{equation*}
    \Vt(\omega) = \inf\left\{S_n \At\circ \psit^{-n}(\zeta, [\omega, \zeta]) + \Delta^s(\omega, [\omega, \zeta]) \: \mid \: n \geq 0, \zeta \in \Omega, \zeta_0 = \omega_0\right\}.
\end{equation*}
\end{definition}

The definition of $m(A)$ guarantees that $\Vt$ is finite. The next proposition will allow us to complete the goal of this section.

\begin{proposition}\label{prop:VisHolder}
    \begin{enumerate}
        \item[(i)] When $\Vt$ is restricted to $\Omega \times \Sigma$, $\At \geq \Vt \circ \psit - \Vt$.
        \item[(ii)] There exists $\mathcal{V}: \Sigma \to \R$ such that $\Vt(\omega) = \mathcal{V} \circ \pi(\omega)$ for all $\omega \in \Omega$.
        \item[(iii)] $\mathcal{V}: \Sigma \to \R$ is globally H\"older. 
    \end{enumerate}
\end{proposition}

Before proving Proposition \ref{prop:VisHolder}, we see that Proposition \ref{prop:DiscretizedSubactionV} then follows easily from it.

\begin{proof}[Proof of Proposition \ref{prop:DiscretizedSubactionV}] We get the Holder map $\V$ from Proposition \ref{prop:VisHolder}, which satisfies $\Vt\circ \tilde{\theta}(\gamma) = \Vt(\theta^\gamma) = \V \circ \pi(\theta^\gamma) = \V(\gamma)$, so that $\Vt\circ \tilde{\theta} = \V$. 
Since we also have $\tilde{\mathcal{A}} \circ \tilde{\theta} = \mathcal{A}$ and $ \tilde{\theta}\circ \psi = \psit \circ \tilde{\theta}$, we obtain
\begin{align*}
    \mathcal{A}(\gamma) = \tilde{\mathcal{A}} \circ \tilde{\theta}(\gamma) & \geq \tilde{\mathcal{V}} \circ \psit(\tilde{\theta}(\gamma)) - \tilde{\mathcal{V}}(\tilde{\theta}(\gamma)) \\
    & = \tilde{\mathcal{V}}\circ \tilde{\theta}\circ \psi(\gamma) - \mathcal{V}(\gamma) \\
    & = \mathcal{V}(\psi(\gamma)) - \mathcal{V}(\gamma).
\end{align*}
\end{proof}

\begin{proof}[Proof of Proposition \ref{prop:VisHolder}]
    Fix $\omega \in \Omega$ and let $\zeta \in \Omega$ such that $\zeta_0 = \omega_0$. Let 
    $$\xi = (\ldots, \zeta_{-2}, \zeta_{-1} \mid \omega_0, \omega_1, \omega_2, \ldots).$$
    Then $\pi(\xi) = [\omega, \zeta]$ and since $\At(-, \pi(\zeta))$ depends only on the forward encoding of $-$, $$\At(\xi, \pi(\xi)) = \At(\omega, \pi(\xi)).$$
    Then
    \begin{align*}
        & S_n\At\circ \psit^{-n}(\xi, [\omega, \zeta]) + \Delta^s(\omega, [\omega, \zeta])+ \At(\omega, \pi(\omega)) \\
        & =  S_n\At\circ \psit^{-n}(\xi, \pi(\xi)) + \left[\At(\omega, [\omega, \zeta]) - \At(\omega, \pi(\omega)) + \Delta^s(\sigma(\omega), \pi\circ\sigma(\xi))\right]+ \At(\omega, \pi(\omega)) \\
        & = S_n\At\circ \psit^{-n}(\xi, \pi(\xi)) + \At(\xi, \pi(\xi)) + \Delta^s(\sigma(\omega), \pi\circ\sigma(\xi))\\
        & = S_{n+1} \At\circ \psit^{-(n+1)}(\sigma(\xi), \pi\circ\sigma(\xi)) + \Delta^s(\sigma(\omega), \pi\circ\sigma(\xi)).
    \end{align*}
    By the Markov property, $\pi\circ\sigma(\xi) \in W_{loc}^s(\sigma(\omega))$, so we obtain:
    $$\Vt(\omega) + \At(\omega, \pi(\omega)) \geq \Vt\circ \sigma(\omega).$$
  For (ii), let $\omega, \omega' \in \Omega$ with $\omega_0 = \omega_0'$ and $ \pi(\omega) = \pi(\omega')$. For any $\zeta$ with $\omega_0 = \zeta_0 = \omega_0'$, we have $[\omega, \zeta] = [\omega', \zeta]$ so $S_n\At\circ\psit^{-n}(\zeta, [\omega, \zeta]) = S_n\At\circ\psit^{-n}(\zeta, [\omega', \zeta])$. Also, by Lemma \ref{lem:prop12replacement}, 
  $\Delta^s(\omega, [\omega, \zeta]) = \Delta^s(\omega', [\omega', \zeta])$. By definition of $\Vt$, we then have $\Vt(\omega) = \Vt(\omega')$. Therefore, $\Vt(\omega)$ depends on $\omega$ actually only via its unique point $\pi(\omega)$; so we can define a map $\V: \Sigma\to \R$ by $$\V(\pi(\omega)) := \Vt(\omega).$$
For (iii), let $\gamma, \gamma' \in \Sigma$. In the following, the constants $K^*$ and $\Lambda_*^s = -2\tau^* < 0< 2\tau^*=\Lambda_*^u$ are from Lemma \ref{lem:bounds on d(psix, psiy)}, the constants $C, \lambda_*^s := -\tau_*<0< \tau_*=: \lambda_*^u$ are from Lemma \ref{lem:exp contraction}, and the constant $\delta$ is from Lemma \ref{lem:lemma24replacement}.

Since the problems for being H\"older show up at small scales, we assume $d_{GX}(\gamma, \gamma') < \frac{\delta}{K^*}$. Let $N = N(\gamma, \gamma')$ be the unique positive integer satisfying
$$ \frac{\delta}{K^*} e^{-(N+1)(\Lambda_*^u - \lambda_*^s)} \leq d_{GX}(\gamma, \gamma') \leq \frac{\delta}{K^*} e^{-N(\Lambda_*^u - \lambda_*^s)}.$$
By Lemma \ref{lem:exp contraction}, there exist pseudo-orbits $\omega, \omega'$ such that 
$\pi(\omega) = \gamma, \: \pi(\omega') = \gamma', \: $ and their first $N+1$ symbols coincide: $\omega_0= \omega'_0, \ldots, \omega_N = \omega'_N.$

Let $n \geq 0$, $\zeta \in \Omega$ with $\zeta_0 = \omega_0 = \omega'_0$ and let $\eta = [\omega, \zeta], \eta'=[\omega', \zeta]$. Since $\eta$ and $\eta'$ are in the stable sets of $\omega$ and $\omega'$, respectively, we may assert that $\zeta$ is encoded so that $\eta$ and $\eta'$ hit exactly the same $N+1$ sections: $\zeta_l = \omega_l (=\omega_l')$ for all $l=0, \ldots, N$. Since taking an infimum of
\begin{equation}\label{eq:VisHolder}
    \left[S_n \At\circ\psit^{-n}(\zeta, \eta) + \Delta^s(\omega, \eta) \right] - \left[S_n \At\circ\psit^{-n}(\zeta, \eta') + \Delta^s(\omega, \eta') \right]
\end{equation}
over all $n\geq 0$ and all such $\zeta$ gives us $\Vt(\omega) - \Vt(\omega') = \V(\gamma) - \V(\gamma')$, we want to bound the magnitude of the expression in \eqref{eq:VisHolder} above in terms of $d_{GX}(\gamma, \gamma')$. 

We use the following notations:
$$\gamma_n = \psit^n_\omega(\gamma), \:\: \gamma'_n = \psit^n_\omega(\gamma'),\:\:  \eta_n = \psit^n_\omega(\eta), \:\: \eta'_n = \psit^n_\omega(\eta').$$
The definition of $N$ gives us
$$N+1 \geq \frac{1}{\Lambda_*^u - \lambda_*^s} \ln\left(\frac{\delta}{K^* d_{GX}(\gamma, \gamma')} \right)$$
which gives us
\begin{equation}\label{eq:DefofN}
    e^{N\lambda_*^s} \leq K_0(\delta) d_{GX}(\gamma, \gamma')^\beta, \: \: \: \: e^{-N\lambda_*^u} \leq L_0(\delta) d_{GX}(\gamma, \gamma')^{-\beta}
\end{equation}
for some constants $K_0(\delta) := e^{-\lambda_*^s} \left(\frac{K^*}{\delta}\right)^\beta$, $L_0(\delta) := e^{\lambda_*^u} \left(\frac{K^*}{\delta}\right)^{-\beta}$, and $\beta:= -\frac{\lambda_*^s}{\Lambda_*^u - \lambda_*^s} \in (0, 1)$. 
Since $\gamma, \eta \in W^s_{loc(\omega)}$ and $\gamma', \eta' \in W^s_{loc(\omega')}$, assuming diam($\Sigma) \leq 1$, 
\begin{equation}\label{eq:e}
    d_{GX}(\eta_N, \gamma_N) \leq C e^{N\lambda_*^s}, \:\: d_{GX}(\eta'_N, \gamma'_N) \leq C e^{N\lambda_*^s}.
\end{equation}
By Lemma \ref{lem:bounds on d(psix, psiy)} and equation \eqref{eq:DefofN}, 
\begin{equation}\label{eq:f}
    d_{GX}(\gamma_N, \gamma'_N) \leq \delta e^{N\lambda_*^s}.
\end{equation}
We also obtain from equations \eqref{eq:e}, \eqref{eq:f} and \eqref{eq:DefofN} (and assuming $\delta\leq C$)
\begin{equation}\label{eq:h}
    d_{GX}(\gamma_N, \gamma'_N), \:d_{GX}(\eta_N, \eta'_N), \: d_{GX}(\gamma_N, \eta_N), \:d_{GX}(\gamma'_N, \eta'_N) \leq M_0 d_{GX}(\gamma, \gamma')^{\beta}
\end{equation}
where $M_0 = 3CK_0.$

Even though $\A$ is highly discontinuous, for any given $\omega\in \Omega$, the map $\At(\omega, -)$ is H\"older continuous, since $\At$ only considers pairs $(\omega, -)$ such that $-$ can follow the orbits that $\omega$ gives instructions for. Therefore, we have a H\"older constant $Hold_\alpha(\At) = \sup_{\gamma, \gamma' \in \Sigma, \omega \in \Omega} \left\{ \frac{|\At(\omega, \gamma) - \At(\omega, \gamma')|}{d(\gamma, \gamma')^\alpha} \right\}$, where $\alpha$ is the exponent constant from the H\"older conditions of $A$ and $\taut$.

We now split the expression in \eqref{eq:VisHolder} into a sum of five parts:
\begin{align*}
    \eqref{eq:VisHolder}_1 & = S_n \At^{-n}(\zeta, \eta) - S_n \At^{-n}(\zeta, \eta') \\
    \eqref{eq:VisHolder}_2 & = \sum_{k=0}^{N-1} \left[\At \circ \psit^k(\omega, \eta) -  \At \circ \psit^k(\omega', \eta') \right]\\
    \eqref{eq:VisHolder}_3 & = \sum_{k=0}^{N-1} \left[\At \circ \psit^k(\omega', \gamma') -  \At \circ \psit^k(\omega, \gamma) \right]\\
    \eqref{eq:VisHolder}_4 & = \sum_{k\geq N} \left[\At \circ \psit^k(\omega, \eta) -  \At \circ \psit^k(\omega, \gamma) \right]\\
    \eqref{eq:VisHolder}_5 & = \sum_{k\geq N} \left[\At \circ \psit^k(\omega', \eta') -  \At \circ \psit^k(\omega', \gamma') \right].
\end{align*}
The sum of the first two terms can be rewritten over the first $(n+N)$ backward iterates of $\eta_N$ and $\eta'_N$ and can be estimated as follows:
\begin{align*}
    |\eqref{eq:VisHolder}_1 + \eqref{eq:VisHolder}_2| & = \left|S_{n+N} \At\circ \psit^{-(n+N)}(\sigma^N(\zeta), \eta_N)  - S_{n+N} \At\circ \psit^{-(n+N)}(\sigma^N(\zeta), \eta'_N) \right| \\
    &=\left|\sum_{k=1}^{n+N} \At \circ \psit^{k-n-N}(\sigma^N(\zeta), \eta_N)  - \sum_{k=1}^{n+N} \At \circ \psit^{k-n-N}(\sigma^N(\zeta), \eta'_N) \right|  \\
    & = \left|\sum_{k=1}^{n+N} \left(\At \circ \psit^{k-n-N}(\sigma^N(\zeta), \eta_N)  - \At \circ \psit^{k-n-N}(\sigma^N(\zeta), \eta'_N) \right) \right| \\
    & \leq \sum_{k=1}^{n+N} Hold_\alpha(\At) \:d_{GX}\left(\psit^{k-n-N}(\sigma^N(\zeta), \eta_N), \psit^{k-n-N}(\sigma^N(\zeta), \eta'_N\right)^\alpha  \\
    & \leq Hold_\alpha(\At)  \sum_{k=1}^{n+N} \left|Ce^{-k\lambda_*^u} d_{GX}(\eta_N, \eta'_N)\right|^\alpha \text{\: (by Lemma \ref{lem:exp contraction})}\\
    & \leq Hold_\alpha(\At)  \sum_{k=1}^{n+N} \left|Ce^{-k\lambda_*^u} \right|^\alpha d_{GX}(\eta_N, \eta'_N)^\alpha \\
    & \leq Hold_\alpha(\At)  \sum_{k=1}^{n+N} \left|Ce^{-k\lambda_*^u} \right|^\alpha \left( M_0 d_{GX}(\gamma, \gamma'_N)^\beta\right)^\alpha \\
    & \leq Hold_\alpha(\At)  \sum_{k=1}^{\infty} \left|Ce^{-k\lambda_*^u} \right|^\alpha  M_0^{\alpha\beta} d_{GX}(\gamma, \gamma'_N)^{\alpha\beta} \\
    & = K_2(\delta) \:d_{GX}(\gamma, \gamma'_N)^{\alpha\beta}.
\end{align*}
For the third term, we use $\omega_k = \omega'_k$ for $0\leq k \leq N$ so that by Lemma \ref{lem:bounds on d(psix, psiy)} and equation \eqref{eq:DefofN},
\begin{align*}
    d_{GX}(\gamma_k, \gamma'_k) & \leq K^* e^{k \Lambda_*^u} d_{GX}(\gamma, \gamma') \\
    & \leq K^*  e^{k \Lambda_*^u} \frac{\delta}{K^*} e^{-N(\Lambda_*^u - \lambda_*^s)} \\
    & = \delta e^{e^{-(N-k)\Lambda_*^u}} e^{N\lambda_*^s} e^{N\lambda_*^s} \\
    & \leq \delta e^{e^{-(N-k)\Lambda_*^u}} e^{N\lambda_*^s} K_0(\delta) d_{GX}(\gamma, \gamma')^\beta. 
\end{align*}
Then we have an estimate for the third term as follows:
\begin{align*}
    |\eqref{eq:VisHolder}_3| & \leq \sum_{k=0}^{N-1} Hold_\alpha(\At) d_{GX}(\gamma_k, \gamma'_k)^\alpha) \\
    & \leq \sum_{k=0}^{N-1} Hold_\alpha(\At) \left(\delta e^{e^{-(N-k)\Lambda_*^u}} e^{N\lambda_*^s} K_0(\delta) d_{GX}(\gamma, \gamma')^\beta \right)^\alpha \\
    & \leq Hold_\alpha(\At)  (\delta K_0(\delta))^\alpha \sum_{l=1}^\infty e^{-l \Lambda^u_*} d_{GX}(\gamma, \gamma')^{\alpha\beta} \\
    & = K_3(\delta) d_{GX}(\gamma, \gamma')^{\alpha\beta}.
\end{align*}
Now, since $\gamma, \eta \in W^s_{loc}(\omega)$ and  $\gamma', \eta' \in W^s_{loc}(\omega')$, using Lemma \ref{lem:bounds on d(psix, psiy)} and equation \eqref{eq:h}, the fourth term can be estimated as follows:
\begin{align*}
     |\eqref{eq:VisHolder}_4| & \leq  Hold_\alpha(\At) \sum_{k\geq N} d_{GX}(\eta_k, \gamma_k)^\alpha \\
     & = Hold_\alpha(\At) \sum_{k=0}^\infty d_{GX}(\eta_{N+k}, \gamma_{N+k})^\alpha \\
     & \leq  Hold_\alpha(\At) \sum_{k=0}^\infty \left(Ce^{k\lambda_*^s} d_{GX}(\eta_N, \gamma_N)\right)^\alpha \\
     & \leq Hold_\alpha(\At) C^\alpha\sum_{k=0}^{\infty} e^{k \lambda_*^s\alpha} \left(M_0(\delta) d_{GX}(\gamma, \gamma')^\beta\right)^\alpha \\
     & = K_4(\delta) d_{GX}(\gamma, \gamma')^{\alpha\beta}.
\end{align*}
Similarly, using $\gamma', \eta'$ and $\omega'$ in place of $\gamma, \eta$ and $\omega$ above, we obtain an estimate for the fifth term:
\begin{align*}
     |\eqref{eq:VisHolder}_5|  \leq  K_4(\delta) d_{GX}(\gamma, \gamma')^{\alpha\beta}.
\end{align*}
Therefore, 
$$|\V(\gamma) - \V(\gamma')| \leq (K_2 + K_3 + K_4 + K_4) d_{GX}(\gamma, \gamma')^{\alpha \beta}$$
as desired.
\end{proof}
\noindent

\section{Proof of the Main Theorem}\label{sec:flow}

With Proposition \ref{prop:DiscretizedSubactionV}, we have solved a discretized version of the sub-action problem. Recall that we have a collection of sections $\Sigma$, and a globally H\"older function $\mathcal{V}:\Sigma \to \mathbb{R}$ such that
\[\mathcal{A}(\gamma) \geq \mathcal{V}\circ\psi(\gamma) - \mathcal{V}(\gamma) \mbox{ for all } \gamma \in \Sigma\]
where $\psi$ is the first-return map for $\Sigma$.  To prove Theorem \ref{thm:flow subactions} we will prove the following proposition.

\begin{proposition}\label{prop:flow subaction}
There is a collection of Poincar\'e subsections $\Sigma'\subset \Sigma$ (with first-return map $\psi'$ and first-return time $\tau'$) and a function $H':GX \to \mathbb{R}_{\geq 0}$ which is globally H\"older and smooth in the flow direction satisfying the following integrability condition: 
\[ \int_0^{\tau'(\gamma)}(A-m(A))\circ g_t(\gamma) dt - (\mathcal{V}\circ\psi'(\gamma)-\mathcal{V}(\gamma)) = \int_0^{\tau'(\gamma)} H'\circ g_t(\gamma) dt\]
for all $\gamma\in \Sigma'$.
\end{proposition}

Proposition \ref{prop:flow subaction} proves Theorem \ref{thm:flow subactions} as follows.

\begin{proof}[Proof of Theorem \ref{thm:flow subactions}]
Given $\Sigma'$ and $\gamma\in GX$, let $T'(\gamma)=\inf\{t\geq 0: g_{-t}\gamma\in \Sigma'\}.$ Then, given $H'$, define
\[V(\gamma) = \mathcal{V}\circ g_{-T'(\gamma)}(\gamma)+ \int_{-T'(\gamma)}^0 (A-m(A)-H')\circ g_t\gamma dt.\]

First, we note that if $\gamma_1\in \Sigma'$, then if $\gamma_1=\psi'(\gamma_2)$, then $T'(\gamma_1)=0$ and so, using the definition of $V$ and the integrability condition given in Proposition \ref{prop:flow subaction} we get
\[V(\gamma_1) = \mathcal{V}(\gamma_1) = \mathcal{V}(\gamma_2) + \int_0^{\tau'(\gamma_2)} (A-m(A)-H')\circ g_t\gamma_2 dt.\]
With this, we have the following well-defined expression for $V$:
\begin{equation}\label{eqn:V}V(\gamma) = \mathcal{V}\circ g_{-T}\gamma + \int_{-T}^0 (A-m(A)-H')\circ g_t \gamma dt
\end{equation}
for all $\gamma$ and any $T$ such that $g_{-T}(\gamma)\in \Sigma'$ or, indeed, in $\Sigma$. Fix any such $T$ for a given $\gamma$. Then 
\begin{align}
    V(g_t \gamma) &= \mathcal{V}\circ g_{-(T+t)}g_t (\gamma) + \int_{-(T+t)}^0 (A-m(A)-H')\circ g_s g_t(\gamma) ds \nonumber \\
    & = \mathcal{V}\circ g_{-T}(\gamma) + \int_0^{T+t}(A-m(A)-H')\circ g_s g_{-T}(\gamma) ds. \nonumber
\end{align}
Taking $\frac{d}{dt}|_{t=0}$ of this expression, we get
\[\frac{d}{dt}|_{t=0}V(g_t\gamma) = (A-m(A)-H')(\gamma).\]
Re-arranging terms, this gives the second statement of Theorem \ref{thm:flow subactions} with $H'$ as the $H$ required in that theorem. Integrating this equation over any geodesic segment yields the first statement of Theorem \ref{thm:flow subactions}. 

For the regularity statements, we work with Equation \eqref{eqn:V}. $\mathcal{V}$ is H\"older on $\Sigma$. On an open set around any $\gamma\in \Sigma'$, we can choose $T$ as a function of $\gamma$ so that $g_{-T}\gamma$ belongs to a single section in $\Sigma$. Using the fact that for the sections we are using, constructed via the methods of \cite{CLT20}, the return map under the flow to a section is H\"older (\cite[Prop 4.9]{CLT20}), $g_{-T}\gamma$ is H\"older in $\gamma$. Finally, $A$ is H\"older by assumption and $H'$ is H\"older by Proposition \ref{prop:flow subaction} so $V$ is globally H\"older. Along orbits of the flow, $V$ is differentiable, with the expression for the derivative given above.
\end{proof}

%

\subsection{Constructing subsections}

To prove Proposition \ref{prop:flow subaction}, following the ideas of \cite{LT}, we will construct a nested sequence of Poincar\'e sections. The following Lemma is a basic tool necessary in these arguments.

\begin{lemma}\label{lem:subsections construction}[cf. Lemma 15 in \cite{LT}]
Given a collection of sections $\Sigma$ as in Lemma \ref{lemma:SectionsCanBeAsSmallAsDesired}, there exists a collection of subsections $\Sigma'$ such that 
\begin{itemize}
    \item[(i)] $\overline\Sigma_i' \subset \Sigma_i$,
    \item[(ii)] $\{U_i':= g_{(0,\tau^*)}\Sigma_i'\}$ cover $GX$,
    \item[(iii)] For all $i,j$, if $U_i \cap \Sigma_j \neq \emptyset$, then $U_i' \cap \Sigma_j' \neq \emptyset$.
\end{itemize}
\end{lemma}

\begin{proof}
Since it is a Poincare section, $\{U_i:=g_{(0,\tau^*)\Sigma_i}\}$ is a finite open cover of $GX$. Let $\epsilon$ be a Lebesgue covering number for $\{U_i\}$.

Using the fact that $g_t$ is continuous, given $\epsilon$, there exists $\delta(\epsilon)>0$ such that whenever $d_{GX}(\gamma_1,\gamma_2)<\delta$, and $t\in [0,\tau^*]$, $d_{GX}(g_t(\gamma_1), g_t(\gamma_2))<\epsilon.$

For each $i$, let $\Sigma_i':=\{\gamma \in \Sigma_i: d_{GX}(\gamma_1, (\cup_{t\in(-\alpha,\alpha)}g_t\Sigma_i)^c) >\frac{\delta(\epsilon)}{2} \}$. That is, we `shrink' $\Sigma_i$ by $\delta(\epsilon)/2$. Clearly (i) is satisfied. For (ii), let $\gamma_3\in GX$. For some $i$, $B_\epsilon(\gamma_3) \subset U_i$ since $\epsilon$ is a Lebesgue number for the covering $\{U_i\}$. Write $\gamma_3 = g_t \gamma_1$, for $(\gamma_1,t)\in \Sigma_i\times(0,\tau^*)$. Suppose that $\gamma_1\notin \Sigma_i'$. Then there is some $\gamma_2\notin \cup_{t\in(-\alpha,\alpha)}g_t\Sigma_i$ such that $d_{GX}(\gamma_1,\gamma_2)<\frac{\delta}{2}$. By the choice of $\delta$, $d_{GX}(\gamma_3,g_t \gamma_2)<\epsilon$, implying that $B_\epsilon(\gamma_3)\nsubseteq U_i$, a contradiction. This proves (ii).

For (iii), if $U_i\cap \Sigma_j \neq \emptyset$, there exists an $\epsilon_{ij}>0$ and $\gamma_3\in \Sigma_j$ such that $B_{\epsilon_{ij}}(\gamma_3) \subset U_i \cap \cup_{t\in(-\alpha, \alpha)}g_t\Sigma_j$. By the same argument used for (ii), if we form $\Sigma_i'$ and $\Sigma_j'$ by shrinking $\Sigma_i$ and $\Sigma_j$ by less than $\epsilon_{ij}$ and less than $\delta(\epsilon_{ij}/2)$, then $\gamma_3\in U_i' \cap \Sigma_j'$. Therefore, if we shrink all our $\Sigma_i$ by $\min\{ \delta(\epsilon)/2, \delta(\epsilon_{ij})/2, \epsilon_{ij}\}$, (iii) is satisfied along with (ii).
\end{proof}

Below, let $\Sigma'$ be a subsection of $\Sigma$ as in Lemma \ref{lem:subsections construction}. Let $\tau'$ and $\psi'$ be the corresponding first-return time and first-return map. For $\gamma\in \Sigma'$, define
\[\mathcal{H}'(\gamma) : = \int_0^{\tau'(\gamma)}(A-m(A))\circ g_t\gamma dt - (\mathcal{V}\circ \psi'(\gamma)-\mathcal{V}(\gamma)).\]
By Proposition \ref{prop:DiscretizedSubactionV}, $\mathcal{H}'\geq 0.$ To prove Proposition \ref{prop:flow subaction}, we need to extend $\mathcal{H}'$ to $H'$ defined on all of $GX$.

%

\subsection{A smoothing function}\label{subsec:h}

A key element in the proof of Proposition \ref{prop:flow subaction} is the `smoothing function' $h$ provided by Lemma \ref{lemma:globalfunc} below. This function will allow us to take the values of $\mathcal{H}$, currently concentrated on the section $\Sigma$ and smooth them out over orbits of the geodesic flow.

\begin{lemma}\label{lemma:globalfunc} There exists a globally Lipschitz continuous, smooth along orbits, non-negative function $h: GX \rightarrow \mathbb{R}_{\geq 0}$ that is null in a neighborhood of $\bigcup\limits_{j = 1}^n \Sigma_j$ such that for all $\gamma \in GX$: 
$$\int_0^{\tau(\gamma)} h \circ g_t(\gamma)dt \geq C$$
\noindent for some constant $C > 0$. 
\end{lemma} 

\begin{proof}
We prove the lemma assuming Lemma \ref{lemma:localfunc}. Let $N = GX \setminus \bigg(\bigcup\limits_{j = 1}^n g_{(-\delta, \delta)} \Sigma_j\bigg)$, where $\delta << \tau_{\ast}$. From Definition \ref{def:Markov proper family}, we know: 
\begin{align}
    N &= \bigg(\bigcup\limits_{i = 1}^n g_{(-\alpha, 0)}\text{Int}_{g}(\Sigma_i)\bigg) \cap \bigg(\bigcup\limits_{j = 1}^n g_{(-\delta, \delta)} \Sigma_j\bigg)^c \nonumber \\
    &= \bigcup\limits_{i = 1}^n \bigg(g_{(-\alpha, 0)}\text{Int}_{g}(\Sigma_i) \cap \big(\cup_{j = 1}^n g_{(-\delta, \delta)} (\Sigma_j)\big)^c\bigg). \nonumber
\end{align} 

Set $U_i =: g_{(-\alpha, 0)}\text{Int}_{g}(\Sigma_i) \cap \big(\cup_{j = 1}^n g_{(-\delta, \delta)} (\Sigma_j)\big)^c$. By Lemma \ref{lemma:localfunc}, there exists some $h_i: GX \rightarrow \mathbb{R}_{\geq 0}$ whose support contains $U_i$ and is contained in $N$ that is smooth in the flow direction and Lipschitz continuous. Let $h=\sum_{i=1}^n h_i$. It is smooth in the flow direction, Lipschitz continuous, and null on a neighborhood of $\bigcup_{j=1}^n \Sigma_j.$

Note that $\mathcal{U} = \{ U_i\}$ forms an open cover of the compact space $GX$. Let $\rho>0$ be a Lebesgue number for this cover. Then for every $\gamma\in GX$, there is some $i$ such that $B_\rho(\gamma)\subset U_i$. Applying (4) of Lemma \ref{lemma:localfunc} we find that
\[\int_0^{\tau(\gamma)} h \circ g_t(\gamma)dt \geq \int_0^{\tau(\gamma)} h_i \circ g_t(\gamma)dt \geq C(\rho)>0.\]
$C(\rho)$ depends only on $\rho$, and hence only on the geometry of the sections, not on $\gamma$, so it can serve as the constant $C$.
\end{proof}

We now prove the key lemma necessary for the proof of Lemma \ref{lemma:globalfunc}. 

\begin{lemma}\label{lemma:localfunc}

As before, let $U_i = g_{(-\alpha, 0)}\text{Int}_{g}(\Sigma_i) \cap \big(\cup_{j = 1}^n g_{(-\delta, \delta)} (\Sigma_j)\big)^c$. Then for every $U_i$, there exists a non-negative function $h_i: GX \rightarrow \mathbb{R}_{\geq 0}$ such that: 
    \begin{itemize} 
        \item[(1)] $U_i \subseteq \text{supp}(h_i) \subseteq GX \setminus \big(\cup_{j = 1}^n g_{(-\delta, \delta)}\Sigma_j\big)$; 
        \item[(2)] $h_i$ is Lipschitz continuous; 
        \item[(3)] $h_i$ is smooth in the flow direction;
        \item[(4)] For all $\gamma$ such that $B_r(\gamma) \subset U_i=W_i$, $\int_0^{\tau(\gamma)} (h_i \circ g_t)(\gamma)dt \geq C(r)$ for some constant $C(r)$ depending only on $r$.
    \end{itemize} 
\end{lemma}

In order to show Lemma \ref{lemma:localfunc}, we define the following operation, which is classical in real analysis: 

\begin{definition} Consider a metric space $X$ equipped with a flow $\phi_t$. Let $\varphi: X \rightarrow \mathbb{R}$ and $\psi: \mathbb{R} \rightarrow \mathbb{R}$.
Fix some $x \in X$ and consider the map $(\varphi \circ \phi_t)(x): \mathbb{R} \rightarrow \mathbb{R}$ defined by $t \mapsto \varphi(\phi_t(x))$. We define the \textit{convolution} of $(\varphi \circ \phi_t)(x)$ and $\psi$, denoted $(\varphi \circ \phi_t)(x) \ast \psi: \mathbb{R} \rightarrow \mathbb{R}$, to be:
$$\big((\varphi \circ \phi_t)(x) \ast \psi\big)(s) = \int \limits_{\mathbb{R}} \big((\varphi \circ \phi_{s - t})(x)\big)\psi(t)dt.$$
\end{definition}

Like in the classical case, the $\ast$ operation is symmetric, which we will now prove: 

\begin{lemma}[Convolution is symmetric] Given $\varphi$ and $\psi$ as above, we have that $(\varphi \circ \phi_t)(x) \ast \psi = \psi \ast (\varphi \circ \phi_t)(x)$.     
\end{lemma}
\begin{proof} The proof follows from a straightforward $u$-substitution: 
\begin{align*} 
\big((\varphi \circ \phi_t)(x) \ast \psi\big)(s) &= \int \limits_{-\infty}^{\infty} \big((\varphi \circ \phi_{s - t})(x)\big)\psi(t)dt\\
&= -\int\limits_{\infty}^{-\infty} \big((\varphi \circ \phi_{u})(x)\big)\psi(s - u)du &\text{($u$-substitution $u = s - t$)}\\
&= \int_{-\infty}^{\infty} \big((\varphi \circ \phi_{u})(x)\big)\psi(s - u)du \\
&= \big(\psi \ast (\varphi \circ \phi_t)(x) \big)(s). 
\end{align*}
\end{proof}

We now specify our input functions. 

Let $\mathcal{U}$ be a cover for a metric space $X$. Recall that if any subset $A \subset X$ of diameter less than $L > 0$ is contained in some $U \in \mathcal{U}$, then $L$ is a \textit{Lebesgue number} for $\mathcal{U}$. Furthermore, a cover $\mathcal{U}$ has \textit{multiplicity} at most $k \geq 0$ if any $x \in X$ belongs to at most $k$ members of $\mathcal{U}$. With these definitions in mind, we consider the following: 

\begin{proposition}[Proposition 4.1, \cite{dg}]\label{prop:pou} Let $\mathcal{U}$ be a cover of a metric space $X$ with multiplicity at most $k + 1$ (where $k \geq 0$), and Lebesgue number $L > 0$. For $U \in \mathcal{U}$, define:
$$\varphi_U(x) = \frac{d(x, X \setminus U)}{\sum\limits_{V \in \mathcal{U}} d(x, X \setminus V)}.$$

Then $\{\varphi_U\}_{U \in \mathcal{U}}$ is a partition of unity on $X$ subordinated to the cover $\mathcal{U}$. Moreover, each $\varphi_U$ satisfies, for all $x, y \in X$:
$$\lvert \varphi_U(x) - \varphi_U(y) \rvert \leq \frac{2k + 3}{L}d(x, y).$$
Furthermore, the family $(\varphi_U)_{U \in \mathcal{U}}$ satisfies, for all $x, y \in X$: 
$$\sum\limits_{U \in \mathcal{U}}\lvert \varphi_U(x) - \varphi_U(y) \rvert \leq \frac{(2k + 2)(2k + 3)}{L}d(x, y).$$
\end{proposition}

\begin{proof}[Proof of Lemma \ref{lemma:localfunc}.]  Let $\psi: \mathbb{R} \rightarrow \mathbb{R}$ be a bump function supported on the interval $(-\epsilon, \epsilon)$ with the property that $\int_{-\epsilon}^{\epsilon} \psi(\gamma)dx = 1$, where $\epsilon << \delta$. Let $W_i = g_{(-\alpha, 0)} \text{Int} B_i \cap \big(\cup_{j = 1}^n g_{(-\delta - \epsilon, \delta + \epsilon)} B_j\big)^c$. Let us denote $\varphi_{W_i}$ from Proposition \ref{prop:pou} as $\varphi_i$. 

Given $\gamma \in GX$, we define the function: 

$$h_i(\gamma) := \big((\varphi_i \circ g_t)(\gamma) \ast \psi\big)(0) = \int_\mathbb{R} \big(\varphi_i \circ g_{- t}\big)(\gamma)\psi(t)dt = \int_{\mathbb{R}} (\varphi_i \circ g_{t})(\gamma)\psi(- t)dt.$$

Note that the two integrals are equal by symmetry of convolution.

\

\noindent\textbf{(1) Support of $\mathbf{h_i}$.}  First, observe that since supp$(\psi) = (-\epsilon, \epsilon)$, it follows that: 
$$h_i(\gamma) = \int_{\mathbb{R}} (\varphi_i \circ g_{-t})(\gamma)\psi(t)dt = \int_{-\epsilon}^{\epsilon} (\varphi_i \circ g_{t})(\gamma)\psi(-t)dt$$

Let $\gamma \in U_i$. Then there exists some union of nonempty open intervals $\mathcal{U}$ that includes an open interval around $0$ such that $(\varphi_i \circ g_{-t})(\gamma) \neq 0$ for all $t \in \mathcal{U}$. As a result, 
$$h_i(\gamma) = \int_{\mathcal{U} \cap (-\epsilon, \epsilon)} (\varphi_i \circ g_{-t})(\gamma)\psi(t)dt > 0.$$
This shows that $\gamma \in \text{supp}(h_i)$, so $U_i \subseteq \text{supp}(h_i)$. 

If $\gamma \in \bigcup\limits_{j = 1}^{n} g_{(-\delta, \delta)} B_j$, then for $t \in (-\epsilon, \epsilon)$, $g_t(\gamma) \in \bigcup\limits_{j = 1}^{n} g_{(-\delta - \epsilon, \delta + \epsilon)} B_j$, which is not in the support of $\varphi_i$. Thus, 

$$h_i(\gamma) = \int_{-\epsilon}^{\epsilon} (\varphi_i \circ g_{-t})(\gamma)\psi(t)dt = 0.$$ 

\

\noindent \textbf{(2) Lipschitz continuity.}
Let $\gamma_1, \gamma_2 \in GX$. Assuming that $\gamma_1$ and $\gamma_2$ are chosen so that supp$(\psi) \cap \text{supp}\big((\varphi_i \circ g_{-t})(\gamma_1) - (\varphi_i \circ g_{-t})(\gamma_2)\big)$ is nonempty (otherwise, the inequality is trivial), then for some $T > \epsilon$: 

\begin{align*}
    \lvert h_i(\gamma_1) - h_i(\gamma_2) \rvert &= \bigg\lvert \int_{\mathbb{R}} (\varphi_i \circ g_{-t})(\gamma_1)\psi(t) - (\varphi_i \circ g_{-t})(\gamma_2)\psi(t) dt \bigg\rvert\\
    &= \bigg\lvert \int_{\mathbb{R}} (\varphi_i \circ g_{t})(\gamma_1)\psi(-t) - (\varphi_i \circ g_{t})(\gamma_2)\psi(-t) dt \bigg\rvert\\
    &\leq \int_{\mathbb{R}} \big\lvert \big((\varphi_i \circ g_{t})(\gamma_1) - (\varphi_i \circ g_{t})(\gamma_2)\big)\psi(-t)\big\rvert dt \\
    &\leq \int_{\mathbb{R}} \bigg(\frac{2K + 3}{L}\bigg)d_{GX}(g_{t}(\gamma_1), g_{t}(\gamma_2))\lvert \psi(-t) \rvert dt &\text{(Lemma \ref{prop:pou})}\\
    &\leq \bigg(\frac{2K + 3}{L}\bigg)e^{2T}d_{GX}(\gamma_1, \gamma_2) \int_{\mathbb{R}} \lvert \psi(-t)\rvert dt &\text{(Lemma \ref{lem:geo flow is Lipschitz})}\\
    &= \bigg(\frac{2K + 3}{L}\bigg)e^{2T}d_{GX}(\gamma_1, \gamma_2).
\end{align*}

We now explain our choice of $T > 0$. Notice that since we are integrating over the a subset of supp$(\psi) = (-\epsilon, \epsilon)$, it follows that any $T > \epsilon$ will suffice. \\

\noindent \textbf{(3) Smoothness along the flow direction.} First, we show that $h_i$ is smooth along the flow direction. In order to do so, we show the infinite differentiability of the function (for \textit{any} $s \in \mathbb{R}$, not just $s = 0$):

$$\big((\varphi_i \circ g_t)(\gamma) \ast \psi\big)(s) = \int_{\mathbb{R}} (\varphi_i \circ g_{t})(\gamma)\psi(s - t)dt.$$

The proof proceeds similarly to the classical case. Recall that:
\begin{align*} \frac{\partial}{\partial s} \big((\varphi_i \circ g_t)(\gamma) \ast \psi\big)(s) &= \lim\limits_{h \rightarrow 0} \dfrac{\big((\varphi_i \circ g_t) \ast \psi\big)(\gamma, s + h) - \big((\varphi_i \circ g_t) \ast \psi\big)(\gamma, s)}{h} \\
&= \lim\limits_{h \rightarrow 0} \int_{\mathbb{R}} \big(\varphi_i \circ g_t\big)(\gamma)\bigg(\frac{\psi(s + h - t) - \psi(s - t)}{h}\bigg)dt. 
\end{align*} 

The next step is to apply the Dominated Convergence Theorem. Since $\psi$ is smooth with compact support, all the derivatives of $\psi$ are bounded, so we can set $M > 0$ to be some number such that $\frac{d}{dt} \lvert \psi(t) \rvert \leq M$. We claim that $M(\varphi_i \circ g_t)(\gamma)$ is an appropriate dominating function. Note that it suffices to show $(\varphi_i \circ g_t)(\gamma)$ is integrable. Suppose $\gamma \notin W_i$. Then: 

\begin{align*}
\int_{-\infty}^{\infty} \lvert (\varphi_i \circ g_t)(\gamma) \rvert dt &\leq \int_{-\alpha}^{\alpha} \lvert (\varphi_i \circ g_t)(\gamma)  - \underbrace{\varphi_i(\gamma)}_{0} \rvert dt \\
& \leq \int_{-\alpha}^{\alpha} \frac{2K + 3}{L} d_{G\widetilde{X}}\big(g_t(\gamma), \gamma\big)dt & \text{(Prop \ref{prop:pou})} \\
&= \int_{-\alpha}^{\alpha} \frac{(2K + 3)t}{L} dt < \infty.
\end{align*}

Otherwise, if $\gamma \in W_i$, then simply replace $\varphi_i(\gamma)$ with $(\varphi_i \circ g_{-2\alpha})(\gamma) = 0$, for example, so $t$ in the last line will be replaced with $t + 2\alpha$, in which case the integral is still finite. 

This allows us to finish showing the derivative exists:
\begin{align*}
\frac{\partial}{\partial s} \big((\varphi_i \circ g_t)(\gamma) \ast \psi\big)(s) &= \int_{\mathbb{R}} (\varphi_i \circ g_t)(\gamma)\psi'(s - t)dt.
\end{align*}
By induction, we have that: 
$$\frac{\partial^n \big((\varphi_i \circ g_t)(\gamma) \ast \psi\big)}{\partial s}(s) = \int_{\mathbb{R}} (\varphi_i \circ g_t)(\gamma)\frac{\partial^n}{\partial s}\psi(s - t)dt.$$
Thus, for any $s \in \mathbb{R}$, $\big((\varphi_i \circ g_t) \ast \psi\big)(\gamma, s)$ is smooth. \\

\noindent \textbf{(4) Boundedness of the integral.} It remains to bound $\int_0^{\tau(\gamma)} (h_i \circ g_t)(\gamma)dt$ below under the condition that $B_r(\gamma) \subset U_i=W_i$.

Recall that $\mathcal{U} = \bigg\{ g_{(-\alpha, 0)}B_i \setminus \big(\cup_{j = 1}^n g_{(-\delta - \epsilon, \delta + \epsilon)}B_j\big) \bigg\}_{i = 1}^n$. Note that $D_0 := \max_{\gamma\in GX} \sum\limits_{U \in \mathcal{U}} d(\gamma, GX \setminus U)$ exists as the sum is finite and the diameter of $GX$ is finite.

Then we have:

\begin{align*}
    \int_{0}^{\tau(\gamma)} h_i(\gamma)ds &= \int_{0}^{\tau(\gamma)} \int_{\mathbb{R}} \big(\varphi_i \circ g_{-t} \circ g_s\big)(\gamma)\psi(t)dtds\\
    &= \int_{0}^{\tau(\gamma)} \int_{\mathbb{R}} \big(\varphi_i \circ g_{s - t}\big)(\gamma)\psi(t)dtds = \int_{0}^{\tau(\gamma)} \int_{\mathbb{R}} \big(\varphi_i \circ g_{t}\big)(\gamma)\psi(s - t)dtds\\
    &\geq \int_{0}^{\tau_{\ast}} \big(\varphi_i \circ g_{t}\big)(\gamma)\int_{t - \epsilon}^{t + \epsilon} \psi(s - t)dsdt \geq \int_{0}^{\epsilon} (\varphi_i \circ g_t)(\gamma)\int_{t - \epsilon}^{t + \epsilon} \psi(s - t)dsdt \\
    &\geq \int_{0}^{\epsilon} (\varphi_i \circ g_t)(\gamma)\int_{-\epsilon}^{\epsilon} \psi(u)dudt = \int_{0}^{\epsilon} (\varphi_i \circ g_t)(\gamma)dt \\
    &\geq \int_{0}^{\epsilon} \frac{d_{G\gamma}(g_t(\gamma), GX \setminus W_i)}{D_0}dt.
\end{align*} 

To complete the argument we simply need to bound $\int_0^\epsilon d_{GX}(g_t(\gamma), GX\setminus W_i)$ below. Since $B_r(\gamma)\in W_i$, and the geodesic flow is unit speed for $d_{GX}$, 
\[d(g_t \gamma, GX\setminus W_i) \geq r-|t| \mbox{ for } t\in [-r,r].\]

Consider two cases. 

\noindent \textbf{Case I: $r>\epsilon$.}
\[ \int_0^\epsilon d_{GX}(g_t \gamma, GX\setminus W_i) dt \geq \int_0^\epsilon r-|t|dt 
         = r\epsilon - \frac{\epsilon^2}{2} > \frac{\epsilon^2}{2}. 
\]

\noindent \textbf{Case II: $r\leq \epsilon$.}
\[ \int_0^\epsilon d_{GX}(g_t \gamma, GX\setminus W_i) dt \geq \int_0^r r-|t|dt \\
         = \frac{r^2}{2}. 
\]

Letting $C(r) = \min\{ \frac{\epsilon^2}{2}, \frac{r^2}{2}\} >0$, we have the desired result.
\end{proof}

%
\subsection{Inductive extension}

To extend the discretized sub-action, we will need to augment the concept of simple transitions.

\begin{definition}\label{def:multipletransitions}
Let $i, j \in I$. We say $i \Rightarrow j$ is a \textit{multiple transition} if there exist $\gamma \in \Sigma_i$ and $n \geq 1$ such that $\psi^n(\gamma) \in \Sigma_j$ and $\tau_{ij}(\gamma) = \sum_{k=0}^{n-1}\tau \circ \psi^k(\gamma) < \tau^*$. The \textit{rank} of the multiple transition $i \Rightarrow j$ is the largest $n \geq 1$ such that there exists a chain $i=i_0 \to i_1 \to \cdots \to i_n = j$ of simple transitions of length $n$ starting at $i$ and ending at $j$. 
\end{definition}
Lopes and Thieullen demonstrated (see Lemma 17 in \cite{LT}) that if we can make the diameter $\alpha$ of the subsections $\Sigma \subset \Sigmat$ sufficiently small–—as permitted by Lemma \ref{lemma:SectionsCanBeAsSmallAsDesired}—–then the rank of any multiple transition is bounded above by $\frac{2\tau^*}{\tau_*}$. This allows us to sensibly define $N$ as the maximum rank of any multiple transition for $\{\Sigma_i\}.$

To extend $\mathcal{H}$ to the function $H'$ specified in Proposition \ref{prop:flow subaction}, we follow the inductive scheme of \cite{LT}. The reason for this inductive argument is ensuring the regularity of $H'$. Extending $\mathcal{H}$ to a flow box based on an individual $\Sigma_i$ is straightforward (see Lemma \ref{lem:Hij0} below). But when two of these flow boxes overlap, ensuring regularity of the resulting extension requires care.

We begin by defining sets used in the construction. Let $i\Rightarrow j$ be a rank $n$ transition.
\begin{itemize}
    \item Using Lemma           \ref{lem:subsections construction} let $\{ \Sigma_i^k\}_i$ for $k=0, \ldots, N$ be a sequence of Poincare sections such that for all $i$, $\Sigma_i^0 = \Sigma_i$, and $\Sigma_i^{k+1}\subset \bar \Sigma_i^{k+1} \subset \Sigma_i^k$ for $k=0, \ldots N-1$.
    
    \item Let $\Sigma_{ij}^k = \{ \gamma\in \Sigma_i^k : \psi_{ij}(\gamma) \in \Sigma_j^k\}.$
    
    \item The \emph{flow box of rank $n$ and size $k$} for $i\Rightarrow j$ is 
    \[B^k_{ij} = \{ g_t\gamma : \gamma\in \Sigma_{ij}^k, 0\leq t\leq \tau_{ij}(\gamma)\}.\] 
    For any $A\subset \Sigma_{ij}^k$,
    \[\{ g_t\gamma : \gamma\in A, 0\leq t \leq \tau_{ij}(\gamma)\}\]
    is a \emph{partial flow box of rank $n$}.

    \item Let \[\Sigma_{ijk}^n = \{\gamma\in \Sigma_{ij}^n : g_t\gamma \in \Sigma_k^n \mbox{ for some } 0\leq t \leq \tau_{ij}(\gamma) \}\] and
    \[U_{ijk}^n = \{ g_t\gamma : \gamma\in \Sigma_{ijk}^n, 0\leq t \leq \tau_{ij}(n)\}.\]
    These are the points in $\Sigma_{ij}^k$ (and corresponding partial flow box) which hit $\Sigma_k^n$ between $\Sigma_i$ and $\Sigma_j$. Note that 
    \begin{align}
        U_{ijk}^n =& \{g_t\gamma: \gamma\in \Sigma_{ijk}^n, 0\leq t \leq \tau_{ik}(\gamma) \} \nonumber \\
                & \cup \{g_t\gamma: \gamma\in \Sigma_{ijk}^n, \tau_{ik}(\gamma)\leq t \leq \tau_{ij}(\gamma) \} \nonumber
    \end{align}
    is the union of two partial flow boxes each of rank less than $n$.

    \item Let 
    \[\mathcal{U}^n = \bigcup_{rank(i\Rightarrow j)\leq n} B_{ij}^n\]
    be the union of all flow boxes of size $n$ and rank $\leq n$. Note that $\mathcal{U}^N = GX.$
\end{itemize}

\begin{figure}[h!]
\centerline{\includegraphics[width=0.7\textwidth]{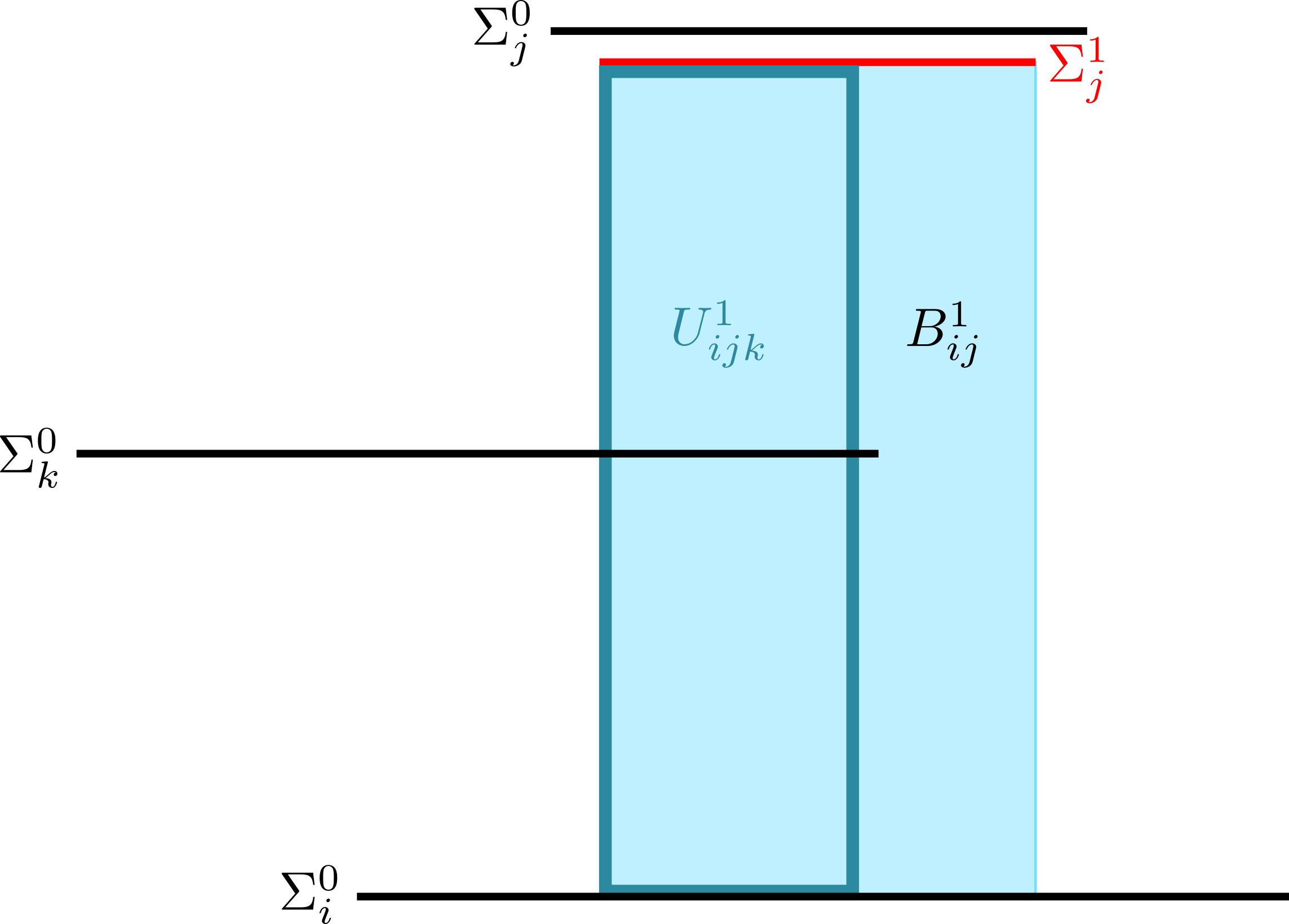}}
	\caption{A illustration of the flow boxes $U_{ijk}^1$ and $B_{ij}^1$.}
	\label{fig:flowboxes}
\end{figure}

We now inductively build $H^n$ on $\mathcal{U}^n$; $H^N$ will be the $H'$ asked for in Proposition \ref{prop:flow subaction}. For this construction we need the following:

\begin{itemize}
    \item For $\gamma\in \Sigma^0_{ij}$, let
\[\mathcal{H}_{ij}(\gamma) := \int_0^{\tau_{ij}^0(\gamma)}(A-m(A))\circ g_t \gamma dt - (\mathcal{V}\circ \psi^0_{ij}(\gamma)-\mathcal{V}(\gamma)).\]
Since $\mathcal{H}_{ij}$ is a sum of $\mathcal{H}\circ \psi^k$, it is nonnegative.
    \item Let $i \Rightarrow j$ be a multiple transition of any rank. On the flow box $B_{ij}^0$, define $H_{ij}^0$ by
\[H_{ij}^0 \circ g_t (\gamma) = \mathcal{H}_{ij}(\gamma) \frac{h\circ g_t (\gamma)}{\int_0^{\tau_{ij}(\gamma)} h\circ g_t(\gamma) dt}\]
for all $\gamma\in \Sigma_{ij}^0$, where $h$ is the function provided by Lemma \ref{lemma:globalfunc}. 
\end{itemize}

\begin{lemma}\label{lem:Hij0}
$H_{ij}^0$ is Lipschitz, smooth in the flow direction, null in a neighborhood of $\bigcup_{j=1}^n \Sigma_j$ and satisfies
\[\mathcal{H}_{ij}(\gamma) = \int_0^{\tau_{ij}(\gamma)} H^0_{ij} \circ g_t (\gamma) dt\]
for all $\gamma\in \Sigma_{ij}^0$.
\end{lemma}

\begin{proof}
The proof follows from the properties of $h$ provided by Lemma \ref{lemma:globalfunc}, the fact that $\tau_{ij}$ is Lipschitz, and straightforward direct computation. 
\end{proof}

The core idea of the construction of $H'$ is in the definition of $H_{ij}^0$.
However, $\{H^0_{ij}\}$ do not jointly define a well-defined function, as there is no reason for them to agree on the overlaps of the $\{B^0_{ij}\}$. Even if this issue is fixed, at the transitions between flow boxes there is no reason for a function patched together from $\{H_{ij}^0\}$ to satisfy the necessary regularity conditions. The inductive construction of \cite{LT} is designed to fix these issues.

\begin{definition}
Let $H^1: \bigcup_{rank(i\to j)=1} B^1_{ij} \to \mathbb{R}$ by $H^1|_{B^1_{ij}}=H^0_{ij}.$
\end{definition}

The only overlaps between the sets $\{B^1_{ij}: rank(i\to j)=1\}$ occur on $\bigcup_i\Sigma_i$. Since $h$ is zero in a neighborhood of this set, $H^1$ is well-defined. It satisfies the integrability condition and regularity conditions thanks to Lemma \ref{lem:Hij0}.

Now suppose that $H^n: \mathcal{U}^n \to \mathbb{R}$ satisfying the integrability and regularity conditions has been defined. Write $H_{ij}^n$ for the restriction of $H^n$ to $B_{ij}^n$.

To define $H^{n+1}$, suppose $rank(i\Rightarrow j)=n+1$. For all $k$ such that $\Sigma_{ijk}^n \neq \emptyset$, as noted above, $U_{ijk}^n$ is a union of two partial flow boxes, each of rank $\leq n$. Hence, $H^n$ is already defined on these partial flow boxes. Therefore, on the partial flow box $\{g_t \gamma: \gamma\in \bigcup_k \Sigma_{ijk}^n, 0\leq t \leq \tau_{ij}(\gamma)\}$, $H^n$ provides a well-defined function $H_{ij}^n$ satisfying the integrability and regularity conditions.

On $\{g_t \gamma: \gamma\in \Sigma_{ij}^n \setminus \bigcup_k \Sigma_{ijk}^n, 0\leq t \leq \tau_{ij}(\gamma)\}$, we have only $H_{ij}^0$. We glue these two functions together with a partition of unity. Let $p,q:\Sigma_i^n \to [0,1]$ be Lipschitz functions such that $p+q=1$ on $\overline{\Sigma}_{ij}^{n+1}$ and so that 
\[\textrm{supp}(p) \subset \bigcup_k \Sigma_{ijk}^{n+1}\]
\[\textrm{supp}(q) \subset  \Sigma_{ij}^n \setminus \bigcup_k \Sigma_{ijk}^{n+1}.\]
Define $H_{ij}^{n+1}$ on $B_{ij}^{n+1}$ by
\[H_{ij}^{n+1}(g_t\gamma) = p(\gamma) H_{ij}^n(g_t\gamma) + q(\gamma)H_{ij}^0(g_t\gamma).\]
Then $H^{n+1}:\mathcal{U}^{n+1}\to \mathbb{R}$ is well-defined by setting $H^{n+1}|_{B_{ij}^{n+1}} = H_{ij}^{n+1}$. It is straightforward to check that it satisfies the integrability condition given that $H_{ij}^n$, $H_{ij}^0$ do and that $p+q=1$ on $\Sigma_{ij}^{n+1}$. The regularity follows from the regularity of $H_{ij}^n$ and $H_{ij}^0$ and the fact that $p$ and $q$ are Lipschitz. This finishes the proof of Proposition \ref{prop:flow subaction}.

\section{Volume Rigidity}\label{sec:MLS}

In the following section, we will show a more general version of the \hyperref[cor:volumerigiditySA]{Volume Rigidity Corollary} stated in the introduction. 

\subsection{Surface Amalgams}
\label{sec:surfaceamalgams}

The following definition is adapted from Definition 2.3 of \cite{LF07}, where they are called \textit{two-dimensional P-manifolds}. 

\begin{definition}[Negatively curved surface amalgams]\label{pmnfld} A compact metric space $X$ is a \textit{negatively curved surface amalgam} if there exists a closed subset $Y \subset X$ (the \textit{gluing curves} of $X$) that satisfies the following: 

\begin{enumerate}
    \item Each connected component of $Y$ is homeomorphic to $S^1$;
    \item The closure of each connected component of $X - Y$ is homeomorphic to a compact surface with boundary endowed with a negatively curved (Riemannian) metric, and the homeomorphism takes the component of $X - Y$ to the interior of a surface with boundary. We will call each $\overline{X - Y}$ a \textit{chamber} in $X$; 
    \item There exists a negatively-curved metric on each chamber which coincides with the original metric. 
\end{enumerate}
\end{definition}

If $Y$ forms a totally geodesic subspace of $X$ consisting of disjoint simple closed curves, we say that $X$ is \textit{simple}. If each connected component of $Y$ (gluing curve) is attached to at least three distinct boundary components of chambers, then we say $X$ is \textit{thick}. Like Lafont in \cite{LF07}, we will only be considering simple, thick surface amalgams, as doing so ensures the surface amalgam is locally CAT($-1$). 

\begin{figure}
    \centering
    \includegraphics[width=\textwidth]{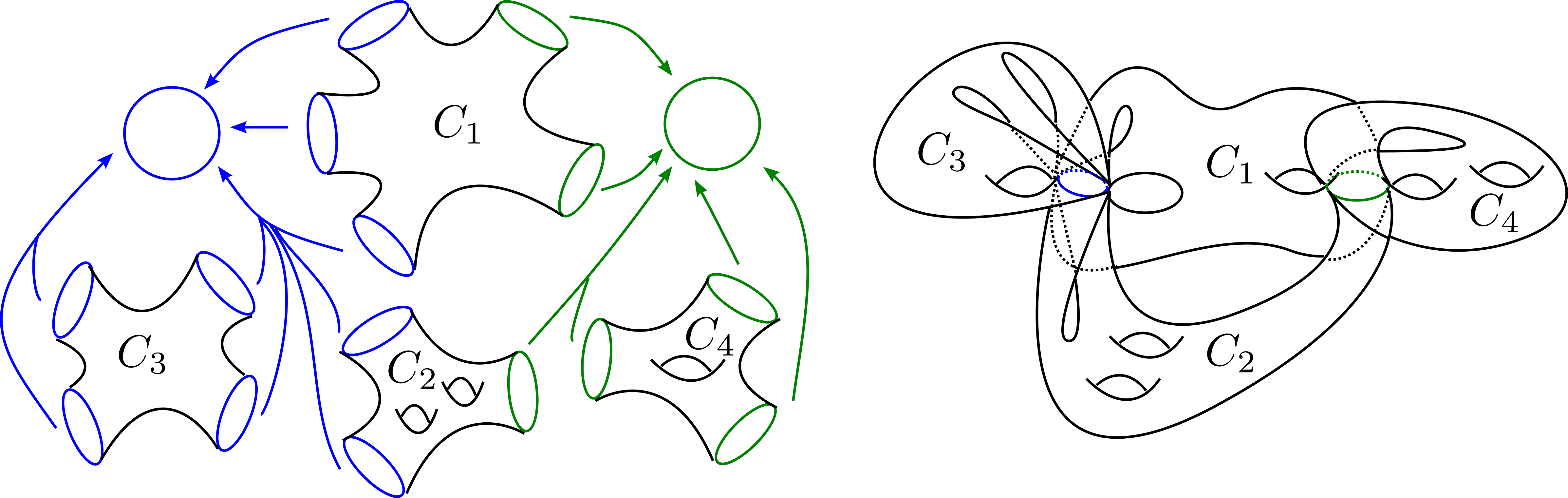}
    \caption{An example of a simple, thick surface amalgam with four chambers.}
    \label{fig:pmnfld}
\end{figure}

We will equip a simple, thick negatively-curved surface amalgam $X$ with a metric in a class we denote as $\mathcal{M}_{\leq}$, following the notation from \cite{cl}. Roughly speaking, metrics in $\mathcal{M}_{\leq}$ are piecewise Riemannian metrics with an additional condition that limits pathological behavior around the gluing curves of $X$. More precisely, we say $g \in \mathcal{M}_{\leq}$ if $g$ satisfies the following properties: 

\begin{enumerate}
    \item Each chamber of $C \subset X$ is equipped with a negatively-curved Riemannian metric with sectional curvature bounded above by $-1$ so that $C$ has geodesic boundary components; 
    \item The restrictions of $g$ to the chambers of $X$ are ``compatible" in the sense that if two boundary components $b_1$ and $b_2$ of two (possibly the same) chambers $C_1$ and $C_2$ are both attached to a gluing curve $\gamma \subset X$, then the gluing maps $b_1 \hookrightarrow \gamma$ and $b_2 \hookrightarrow \gamma$ are isometries (in particular, we do not allow circle maps of degree two;
    \item For any two boundary components $b_1 \in C_1$ and $b_2 \in C_2$, the restriction of $g$ to $N_{b_1} \bigcup\limits_{b_1 \sim b_2} N_{b_2}$ is a negatively curved smooth Riemannian metric with sectional curvature bounded above by $-1$, where $N_{b_1}$ and $N_{b_2}$ are $\epsilon$-neighborhoods around $b_1$ and $b_2$ respectively for some $\epsilon > 0$.
    
\end{enumerate}

We impose the third condition to ensure that we can exploit previous marked length spectrum rigidity results for surfaces which in particular require Riemannian negatively curved metrics with at most a finite number of cone singularities. We now discuss some properties of $\mathcal{M}_{\leq}$ that will be useful in the proof of Theorem \ref{thm:volumerigidity}. 

\begin{remark}\label{cat-1}
If $(X, g)$ is a negatively-curved surface amalgam where $g \in \mathcal{M}_{\leq}$, then $(X, g)$ is locally CAT(-1).  
\end{remark}

Indeed, suppose $X$ is equipped with a metric $g \in \mathcal{M}_{\leq}$ and $C \subset X$ is a chamber in $X$. Recall a generalization of the Cartan-Hadamard Theorem which states that a smooth Riemannian manifold $M$ has sectional curvature $\leq \kappa$ if and only if $M$ is locally CAT($\kappa$) (see \cite{bh} Theorem 1A.6). As a result, the restriction of $g$ to $C$ is locally CAT(-1) since $C$ is a endowed with a negatively curved metric with sectional curvature bounded above by $-1$.  If $\kappa \in \mathbb{R}$ and $X_1$ and $X_2$ are locally CAT($\kappa$) spaces glued isometrically along a convex, complete metric subspace $A \subset X_1 \cap X_2$, then $X_1 \sqcup_A X_2$ is locally CAT($\kappa$) (see Theorem 2.11.1 in \cite{bh}). As a result, a negatively curved surface amalgam $(X, g)$ with locally CAT(-1) chambers will also be locally CAT(-1), as claimed.

%
\subsection{Volume Rigidity}

 Let $\mathscr{G}\widetilde X$ be the set of unparametrized, un-oriented geodesics in $\widetilde X$. A \textit{geodesic current} on $(X, g)$ is a $\pi_1$-invariant Radon measure on $\mathscr{G}\widetilde{X}$. 
 
 There are two especially important examples of geodesic currents. We will assume, for simplicity, that $(X, g)$ is locally CAT($-1$). Furthermore, we assume there is a well-defined notion of \textit{transversality} in $\mathscr{G}\widetilde X$. If $(X, g)$ is a negatively-curved surface, there is a natural such notion since $\partial^\infty\widetilde{X}$ is a circle. For other settings, transversality must be defined with more care; for examples, see Section 6 of \cite{cl} and Section 2.4.1 of \cite{wu}. 

First, given a homotopy class $[\alpha] \in \pi_1(X)$ with geodesic representative $\alpha$, there is a \textit{counting current} associated with $\alpha$ that assigns to each Borel set $E \subset \mathscr{G}\widetilde{X}$ the measure: 
$$\mu_{\alpha}(E) = \lvert E \cap \{\widetilde{\alpha}\} \rvert$$
where $\{\widetilde{\alpha}\}$ denotes the set of lifts of $\alpha$ to $\widetilde{X}$. We follow convention and, with abuse of notation, write $\alpha := \mu_\alpha$, which makes clear the fact that the collection of closed geodesics of $(X, g)$ embeds naturally into the space of geodesic currents.

Second, there is a \textit{Liouville current} which, roughly speaking, captures lengths of closed geodesics. There is no general definition of a Liouville current; rather, the Liouville current should capture important properties, which we state in Assumption \ref{assumption:liouville}. Before stating the assumption, we need to introduce the notion of \textit{intersection numbers}. 

In the case of negatively curved Riemannian surfaces, the geometric intersection number of two closed geodesics (viewed as geodesic currents) extends to a symmetric, bilinear form (see Proposition 4.5 of \cite{bonahon}), the \textit{intersection number} of two geodesic currents. Explicitly, given two geodesic currents $\mu, \nu$, one can define the intersection number of $\mu$ and $\nu$ as:
$$i(\mu, \nu) = \int_{DGX} d\mu \times d\nu$$
where $DGX$ is the set of all pairs of transversally intersecting geodesics in $GX$.

A straightforward computation shows that for two counting currents $[\alpha]$ and $[\beta]$, $i(\alpha, \beta)$ is exactly their geometric intersection number. Metric-independent generalizations of intersection numbers can be found for compact quotients of certain Fuchsian buildings (see \cite{cl}) and surface amalgams (see \cite{wu}). 

We now state a few requirements for the Liouville current necessary for the proof of Theorem \ref{thm:volumerigidity}. 

\begin{assumption}\label{assumption:liouville} Let $(X, g)$ be a locally CAT(-1) metric space with a well-defined notion of transversality, a well-defined intersection number function as well as a \textit{Liouville current} $\Lambda_g$, e.g. a geodesic current that satisfies the following two properties:

\begin{enumerate}
    \item $i(\alpha, \Lambda_g) = C_\alpha\ell_g(\alpha) = C_\alpha\mathscr{L}_{g}([\alpha])$; 
    \item $i(\Lambda_g, \Lambda_g) = K\pi\text{Vol}_g(X)$
\end{enumerate}

\noindent where $C_\alpha$ and $K$ are two positive real constants.
\end{assumption}

\

Assumption \ref{assumption:liouville} holds in the case of negatively curved Riemannian surfaces (see \cite{otal}), and for negatively curved surfaces with large angle cone points as in \cite{HP}. In this case (and with Otal's choice of scale factor) $C_\alpha = 1$ for all $\alpha$ and $K = \frac{1}{2}$. A Liouville measure, a flow-invariant Radon measure on \textit{oriented} unit-speed geodesics, is defined for nonpositively curved orbihedra in \cite{bb} which is used in \cite{cl} and \cite{wu} for quotients of Fuchsian buildings and surface amalgams respectively. Both sets of authors choose to scale the Liouville current, viewed as a measure on \textit{unoriented} unit-speed geodesics, to match the scale factor of the Liouville measure defined in \cite{bb}. However, if one scales the Liouville measure from \cite{bb} with a factor of $\frac{1}{2}$ to match the scale factor in \cite{otal}, one can see that for the cases of surface amalgams and quotients of Fuchsian buildings, $C_\alpha = 2$ and $K = 1$, given that $\alpha$ is \textit{not} a branching geodesic or gluing curve. (We remark that the original factor $K = 4$ from \cite{cl} should be $K = 2$ and has been calculated without the scale factor of $\frac{1}{2}$). 
So Assumption \ref{assumption:liouville} is satisfied.


We also require the intersection pairing to satisfy a specific form of continuity:

\begin{assumption}[Weak continuity of the intersection number function] \label{assumption:liouville3} Given $g$ and $g'$, two metrics on $X$, and a sequence of currents $(\mu_k)_{k \in \mathbb{N}}$ on $X$ such that $\mu_k \rightarrow \Lambda_g$ with respect to the weak-$\ast$ topology, $i(\mu_k, \Lambda_{g'}) 
\rightarrow i(\Lambda_g, \Lambda_{g'})$.  
\end{assumption}

In the case of closed nonpositively-curved surfaces (with or without cone points), the intersection number given in \cite{bonahon} is continuous everywhere and thus automatically satisfies Assumption \ref{assumption:liouville3}. In the cases of compact quotients of certain Fuchsian buildings and surface amalgams, the intersection number functions are \textit{not} continuous everywhere (see Lemma 10.2 of \cite{cl}). However, by Proposition 10.1 in \cite{cl}, they satisfy Assumption \ref{assumption:liouville3}; the proof follows verbatim in the case of surface amalgams. 

Finally, in order to prove Theorem \ref{thm:volumerigidity}, we need one final assumption about the Liouville current:

\begin{assumption}\label{assumption:liouville2} $(X, g)$ is equipped with a Liouville current that is realized as the weak-$\ast$ limit of scalar multiples of counting currents $a_n\alpha_n$ with $C_{\alpha_n}=2K$.
\end{assumption}

In the cases of negatively-curved closed surfaces (with or without cone points), quotients of Fuchsian buildings from \cite{cl}, and surface amalgams from \cite{wu}, the following property is satisfied: scalar multiples of counting currents are dense in the space of all geodesic currents equipped with the weak-$\ast$ topology. One can see this is in a number of different ways. 

Bonahon proves this statement in the setting of Gromov hyperbolic spaces in Theorem 7 of \cite{bonahon91}. While he works with Cayley graphs of Gromov hyperbolic groups, his argument mainly relies on the existence of a free, cocompact, properly discontinuous, isometric action on the space as well as properties of quasigeodesics in CAT(-1) and Gromov hyperbolic spaces. 

Alternately, as noted in \cite[Proposition 2]{bonahon_geometry}, this statement follows from \cite[Theorem 1]{sigmund} or \cite[Theorem 1]{sigmundflow}, which only use the weak specification property of the geodesic flow.  Weak specification holds for geodesic flow on any compact, locally CAT(-1) space \cite[Theorem A]{clt}. Sigmund's density proof uses weak-specification to approximate the Liouville measure with measures supported on closed geodesics. In Fuchsian buildings and surface amalgams, the only $\alpha$ for which $C_\alpha \neq 2K$ are the branching geodesics or gluing curves. These are proper, closed subsets in $GX$ and it is easy to see from the definition of the Liouville measure developed by \cite{bb, cl, wu} that such geodesics have zero Liouville measure. Therefore, any such geodesics can be omitted from the approximating sequence $(a_n\alpha_n)\to \Lambda_g$ obtained via Sigmund's argument. Therefore the examples we are interested in satisfy Assumption \ref{assumption:liouville2}.

Recall that a class of CAT($-1$) metrics $\mathcal{M}$ is \textit{marked length spectrum rigid} if whenever $(X, g_1)$ and $(X, g_2)$ (where $g_1, g_2 \in \mathcal{M}$) have the same marked length spectrum, $(X, g_1)$ and $(X, g_2)$ are isometric. There are some well-documented classes of marked length spectrum rigid metrics. Negatively curved Riemannian metrics on surfaces are known to be marked length spectrum rigid due to $\cite{croke}$ and $\cite{otal}$. Due to \cite{cl}, certain classes of piecewise negatively curved Riemannian metrics, including piecewise hyperbolic metrics, on compact quotients of Fuchsian buildings are marked length spectrum rigid. Piecewise negatively curved Riemannian metrics satisfying certain smoothness conditions on simple, thick surface amalgams are also marked length spectrum rigid due to \cite{wu}. Marked length spectrum rigidity is a key tool in proving volume rigidity; thus, we will state it as an additional assumption. 

Recall also from Section \ref{sec:intro} that given a metric $g$ on a metric space $(X, g)$, $\mathscr{L}_g$ denotes the marked length length spectrum of $(X, g)$. We are now ready to state and prove our volume rigidity result, which is a generalization of the \hyperref[cor:volumerigiditySA]{Volume Rigidity Corollary} stated in the introduction: 

\begin{theorem}[Volume Rigidity]\label{thm:volumerigidity}
Let $(X, g_0)$ and $(X, g_1)$ be two locally CAT(-1) spaces satisfying Assumptions \ref{assumption:liouville}, \ref{assumption:liouville3} and \ref{assumption:liouville2}. Furthermore, suppose $g_0$ and $g_1$ belong to a class of marked length spectrum rigid metrics. Let $\mathscr{L}_{g_0} \leq \mathscr{L}_{g_1}$. Then Vol$_{g_0}(X) \leq \text{Vol}_{g_1}(X)$. Furthermore, if Vol$_{g_0} = \text{Vol}_{g_1}$, then $(X, g_0)$ and $(X, g_1)$ are isometric. 
\end{theorem} 

\begin{proof}
There exists a $\pi_1(X)$-equivariant homeomorphism $\Phi: (G\Xt, \gt_0) \to (G\Xt, \gt_1)$. Let $\Phi^*\Lambda_{g_1}$ be the pullback of $\Lambda_{g_1}$ under $\Phi$ so that $\Phi^*\Lambda_{g_1}$ is a geodesic current of $(X, g_0)$. (Technically, we are using the homeomorphism induced by $\Phi$ on unoriented, unparameterized geodesics.) By Assumption \ref{assumption:liouville}(1) and the hypothesis, for any $[\alpha] \in \pi_1(X)$, 
    \begin{align*}
                i(\alpha, \Lambda_{g_0}) =  C_\alpha\mathscr{L}_{g_0}([\alpha])  & \leq  C_\alpha \mathscr{L}_{g_1}([\alpha]) \\
                & = i(\alpha, \Lambda_{g_1})= i( \Phi^*\alpha, \Phi^* \Lambda_{g_1}) = i(\alpha, \Phi^* \Lambda_{g_1}).
    \end{align*}
Using the fact that the Liouville current is realized as a limit of counting currents (Assumption \ref{assumption:liouville2}) and the continuity of $i(-, -)$ at the Liouville current (Assumption \ref{assumption:liouville3}), 
$$i(\Lambda_{g_0}, \Lambda_{g_0}) \leq i(\Lambda_{g_0}, \Phi^* \Lambda_{g_1}).$$

Finally, using the symmetry of the intersection form, another application of Assumption \ref{assumption:liouville2}, and  Assumption \ref{assumption:liouville}(2), 
$$K\pi \text{Vol}_{g_0} = i(\Lambda_{g_0}, \Lambda_{g_0}) \leq i(\Lambda_{g_0}, \Phi^* \Lambda_{g_1}) \leq i(\Phi^* \Lambda_{g_1}, \Phi^*\Lambda_{g_1}) = K\pi \text{Vol}_{g_1}$$
and we get the desired volume inequality. Note that the condition on $C$ and $K$ in Assumption \ref{assumption:liouville2} is not necessary for this part of the proof.

Now suppose there is a volume equality: $\text{Vol}_{g_0} = \text{Vol}_{g_1}$. Then, in particular, we have the following equality:
\begin{equation}\label{eq:volequality}
    i(\Lambda_{g_0}, \Phi^* \Lambda_{g_1}) = K\pi \text{Vol}_{g_0}.
\end{equation}
Let $g^0_t$ and  $g^1_t$ be the geodesic flows on $(GX, g_0)$ and $(GX, g_1)$, respectively. Let $F: (GX, g_0) \to (GX, g_1)$ be an orbit equivalence of these flows. Then there exists a reparameterization map $T: (GX, g_0) \times \R \to \R$ such that we have
$$F(g^0_t(\gamma)) = g^1_{T(\gamma, t)}(F(\gamma)).$$  By the general stability theory (see \cite{KH95}, Chapter 19), the orbit equivalence $F$ (and consequently $T$) can be chosen to be H\"older continuous. The key requirements for this are the existence of stable and unstable foliations  (Theorem \ref{theorem:GeodesicFlowMetricAnosov}) for the flow and the Lipschitz continuity of the geodesic flows (Lemma \ref{lem:geo flow is Lipschitz}).

Since $T(\gamma, s+1) = T(\gamma, s) + T(g^0_s(\gamma), 1),$ for any $[\alpha] \in \pi_1(S)$ if we denote $L:= \mathscr{L}_{g_0}([\alpha])$, then
\begin{align*}
    \int_0^{L} T(g^0_s(\alpha), 1) \; ds & = \int_0^{L} T(\alpha, s+1) - T(\alpha, s) \; ds\\
    & = \int_0^1 T(\alpha, L+s) - T(\alpha, s) \; ds \\
    & = T(g^0_L(\alpha), L) \\
    & = T(\alpha, L) = \mathscr{L}_{g_1}([\alpha])\\
    & = \frac{1}{C_\alpha} i(\alpha, \Phi^*\Lambda_{g_1}).
\end{align*}
Therefore, for any $[\alpha] \in \pi_1(X)$, 
\begin{equation}\label{eqn:pair with pullback}
    i(\alpha, \Phi^*\Lambda_{g_1}) = C_\alpha \int_0^{\mathscr{L}_{g_0}([\alpha])} T(g_s^0(\alpha), 1)\; ds.
\end{equation}

We now want to interpret the right-hand side of the above equation as an integral over $GX$ with respect to a geodesic-flow invariant measure. This requires some care.

As before, let $\partial^{(2)}_{\infty}\tilde X = (\partial^\infty\tilde X \times \partial^\infty\tilde X) \setminus \Delta$ be the set of distinct, ordered pairs of points in $\partial^\infty\tilde X$. This set specifies the \emph{oriented}, unparamerized geodesics in $\tilde X$ and there is a natural 2-to-1 map $\Psi:\partial^{(2)}_{\infty}\tilde X \to \mathscr{G}\tilde X$ sending $(\xi,\eta)\mapsto \{\xi,\eta\}$ mapping to the \emph{unoriented}, unparametrized geodesics $\mathscr{G}\tilde X$. For any measure $\mu$ on $\partial^{(2)}_{\infty}\tilde X$, the push-forward $\Psi_*\mu$ is a measure on $\mathscr{G}\tilde X$ defined by $\Psi_*\mu(A) = \mu(\Psi^{-1}A).$ If $\mu$ is $\Gamma$-invariant, so is $\Psi_*\mu$, and hence it is a current. 

Let $\tilde \Lambda_g$ be the measure on $\partial^{(2)}_{\infty}\tilde X$ defined in local coordinates along a geodesic segment exactly as the Liouville current is. For instance, on a surface its expression would be $\frac{1}{2}|\sin\theta| d\theta dx$ where $\theta\in[0,2\pi]$, with a corresponding expressions for a Fuchsian building or surface amalgam as described in \cite{wu} or \cite{cl} (up to the choice of scale factor $\frac{1}{2}$). Note that since this measure is defined on the space of \emph{oriented} geodesics, the angular coordinate must run over $[0,2\pi]$ instead of its domain for $\mathscr{G}\tilde X$, $[0,\pi]$. 

Since $\Psi$ is 2-to-1,
\[\Psi_*\tilde \Lambda_g = \sin\theta d\theta dx = 2\Lambda_g\]
(where now $\theta\in [0,\pi]$). If $\alpha\in \pi_1(X)$, then (again abusing notation a bit) let $\alpha$ be the measure on $\partial^{(2)}_{\infty}\tilde X$ given by the counting measure on \emph{oriented} lifts of the geodesic representative of $\alpha$, with the orientation induced by $\alpha$ -- namely oriented from the repelling to the attracting endpoint of (the correct conjugate of) $\alpha$ for its action on $\partial^\infty\tilde X$. Then $\Psi_*\alpha = \alpha$ where on the right-hand side of this expression $\alpha$ is the counting measure on \emph{unoriented} lifts of $\alpha$. This equation holds because each unoriented lift of $\alpha$ has exactly one pre-image under $\Psi$ in the support of $\alpha$ as a measure on $\partial^{(2)}_{\infty}\tilde X$, namely the same lift oriented according to the action of (the correct conjugate of) $\alpha$.

Under Assumption \ref{assumption:liouville2}, we can take a sequence $a_n\alpha_n \to \tilde \Lambda_{g_0}$ in the weak-$*$ topology for measures on $\partial^{(2)}_{\infty}\tilde X$ for $a_n\in\mathbb{R}^+$, $\alpha_n\in\pi_1(X)$, and $C_{\alpha_n}=2K$. Using the notes of the previous paragraph, 
\[\Psi_*(a_n\alpha_n)=a_n\alpha_n \to \Psi_*\tilde \Lambda_{g_0} = 2\Lambda_{g_0}\]
in the weak-$*$ topology for geodesic currents.

Consider the measure $\tilde\Lambda_{g_0}\times dt$ on $\partial^{(2)}_{\infty}\tilde X \times \mathbb{R} \cong G\tilde X$. In \cite[(A3)]{bonahon91}, Bonahon examines this measure in its local coordinates and shows that it is $\frac{1}{2}\lambda_g$ where $\lambda_g$ is the Liouville measure. Returning to Equation \eqref{eqn:pair with pullback}, using continuity of the intersection pairing at a Liouville current, the fact that when $a_n\alpha_n\to\tilde\Lambda_{g_0}$, $a_n\alpha_n\to 2\Lambda_{g_0}$, and the continuity of $T$, we get
\begin{align}
    2i(\Lambda_{g_0},\Phi^*\Lambda_{g_1}) & = \lim_{n\to\infty} C_{\alpha_n}a_n \int_0^{\mathscr{L}([\alpha_n])}T(\alpha_n,1)ds \nonumber \\
    & = 2K \lim_{n\to\infty} \int_{GX\cong (\partial^{(2)}\tilde X\times \mathbb{R})/\Gamma} T(\gamma,1) d(a_n\alpha_n)dt \nonumber \\
    & = 2K \int_{GX} T(\gamma,1) d\tilde \Lambda_{g_0}dt \nonumber \\
    & = \frac{2K}{2}\int_{GX} T(\gamma,1)d\lambda_{g_0}. \nonumber
\end{align}
Therefore,

\begin{equation}\label{eq:EquationT(gamma,1)}
    i(\Lambda_{g_0}, \Phi^*\Lambda_{g_1}) = \frac{K}{2} \int_{GX} T(\gamma, 1) \; d\lambda_{g_0}.
\end{equation}
Using equations \eqref{eq:EquationT(gamma,1)} and \eqref{eq:volequality},
\begin{align}\label{eq:intA=0}
    \int_{GX} \big(T(\gamma, 1) - 1 \big)\; d\lambda_{g_0} & = \frac{2}{K} i(\Lambda_{g_0}, \Phi^*\Lambda_{g_1}) - 2\pi \text{Vol}_{g_0}  \notag \\ 
    & = \frac{2}{K} K\pi \text{Vol}_{g_0} - 2\pi \text{Vol}_{g_0}  = 0.
\end{align}

We now use the \hyperref[thm:flow subactions]{Main Theorem} with $A: GX \to \R$ defined as 
$A(\gamma) := T(\gamma, 1) - 1.$ We first claim that the integral of $A$ along any closed geodesic $\alpha$ is zero. By the \hyperref[thm:flow subactions]{Main Theorem}, we have a sub-action $V: GX \to \R$ that is smooth in the flow direction such that $$A(\gamma) = m + \left( \frac{d}{dt}\right)\Big|_{t=0} V(g^0_t(\gamma)) + H(\gamma).$$ Note that $H$ is a non-negative function. Moreover, by the marked length spectrum inequality assumption, for any closed geodesic $\alpha$, $T(\alpha, 1) \geq 1$, so using Sigmund's Theorem \cite[Theorem 1]{sigmundflow}, we get that the minimal average $m(A)$ of $A$ is also non-negative. Therefore, using the Lie derivative notation $L_X V(\gamma) : = \left( \frac{d}{dt}\right)\Big|_{t=0} V(g^0_t(\gamma))$, we have $A \geq L_X V$. Together with \eqref{eq:intA=0}, this gives 
\begin{align*}
    0 & \leq \int_{GX} \big(A(\gamma) - L_XV(\gamma)\big)\; d\lambda_{g_0} = \int_{GX} A(\gamma) \; d\lambda_{g_0} = 0,
\end{align*}
and hence that $A(\gamma) = L_X V(\gamma)$. For any closed geodesic $\alpha$, 
$$\int_0^{\ell_{g_0}(\alpha)} T(\alpha, 1) - 1 \; dt = \int_0^{\ell_{g_0}(\alpha)} A(g^0_t(\alpha))\; dt = \int_0^{\ell_{g_0}(\alpha)} L_X V(g^0_t(\alpha))\; dt = 0,$$ 
as claimed.
Now, for any free homotopy class $[\alpha] \in \pi_1(X)$, 
$$\mathscr{L}_{g_0}([\alpha]) =  \int_0^{\ell_{g_0}(\alpha)} 1 \; dt = \int_0^{\ell_{g_0}(\alpha)} T(\alpha, 1) \; dt = \mathscr{L}_{g_1}([\alpha]).$$
Finally, by our assumption marked length spectrum rigidity for the class of metrics containing $(X,g_i)$, $g_0$ and $g_1$ are isometric. 
\end{proof} 
\bibliographystyle{amsalpha}
\bibliography{References}

\providecommand{\bysame}{\leavevmode\hbox to3em{\hrulefill}\thinspace}
\providecommand{\MR}{\relax\ifhmode\unskip\space\fi MR }
\providecommand{\MRhref}[2]{%
  \href{http://www.ams.org/mathscinet-getitem?mr=#1}{#2}
}
\providecommand{\href}[2]{#2}
\begin{thebibliography}{CLT20b}

\bibitem[BB95]{bb}
Werner Ballmann and Michael Brin, \emph{Orbihedra of nonpositive curvature},
  Inst. Hautes \'Etudes Sci. Publ. Math. (1995), no.~82, 169--209. \MR{1383216}

\bibitem[BH99]{bh}
Martin~R. Bridson and Andr\'{e} Haefliger, \emph{Metric spaces of non-positive
  curvature}, Grundlehren der mathematischen Wissenschaften [Fundamental
  Principles of Mathematical Sciences], vol. 319, Springer-Verlag, Berlin,
  1999.

\bibitem[Bon86]{bonahon}
Francis Bonahon, \emph{Bouts des vari\'et\'es hyperboliques de dimension
  {$3$}}, Ann. of Math. (2) \textbf{124} (1986), no.~1, 71--158. \MR{847953}

\bibitem[Bon88]{bonahon_geometry}
\bysame, \emph{The geometry of {T}eichm\"uller space via geodesic currents},
  Invent. Math. \textbf{92} (1988), no.~1, 139--162. \MR{931208}

\bibitem[Bon91]{bonahon91}
\bysame, \emph{Geodesic currents on negatively curved groups}, Arboreal group
  theory ({B}erkeley, {CA}, 1988), Math. Sci. Res. Inst. Publ., vol.~19,
  Springer, New York, 1991, pp.~143--168. \MR{1105332}

\bibitem[Bou95]{bourdon}
Marc Bourdon, \emph{Structure conforme au bord et flot g\'{e}od\'{e}sique d'un
  {${\rm CAT}(-1)$}-espace}, Enseign. Math. (2) \textbf{41} (1995), no.~1-2,
  63--102.

\bibitem[Bow73]{Bow73}
Rufus Bowen, \emph{Symbolic dynamics for hyperbolic flows}, Amer. J. Math.
  \textbf{95} (1973), 429--460. \MR{339281}

\bibitem[Bro41]{brooks}
R.~L. Brooks, \emph{On colouring the nodes of a network}, Proc. Cambridge
  Philos. Soc. \textbf{37} (1941), 194--197. \MR{12236}

\bibitem[CD04]{CD04}
Christopher~B. Croke and Nurlan~S. Dairbekov, \emph{Lengths and volumes in
  {R}iemannian manifolds}, Duke Math. J. \textbf{125} (2004), no.~1, 1--14.
  \MR{2097355}

\bibitem[CL19]{cl}
David Constantine and Jean-Fran\c{c}ois Lafont, \emph{Marked length rigidity
  for {F}uchsian buildings}, Ergodic Theory Dynam. Systems \textbf{39} (2019),
  no.~12, 3262--3291. \MR{4027549}

\bibitem[CLT20a]{CLT20}
David Constantine, Jean-Fran\c{c}ois Lafont, and Daniel~J. Thompson,
  \emph{Strong symbolic dynamics for geodesic flows on {${\rm CAT}(-1)$} spaces
  and other metric {A}nosov flows}, J. \'{E}c. polytech. Math. \textbf{7}
  (2020), 201--231. \MR{4054334}

\bibitem[CLT20b]{clt}
\bysame, \emph{The weak specification property for geodesic flows on {${\rm
  CAT}(-1)$} spaces}, Groups Geom. Dyn. \textbf{14} (2020), no.~1, 297--336.

\bibitem[Coo93]{coornaert}
Michel Coornaert, \emph{Mesures de {P}atterson-{S}ullivan sur le bord d'un
  espace hyperbolique au sens de {G}romov}, Pacific J. Math. \textbf{159}
  (1993), no.~2, 241--270. \MR{1214072}

\bibitem[Cro90]{croke}
Christopher~B. Croke, \emph{Rigidity for surfaces of nonpositive curvature},
  Comment. Math. Helv. \textbf{65} (1990), no.~1, 150--169.

\bibitem[DG07]{dg}
Marius Dadarlat and Erik Guentner, \emph{Uniform embeddability of relatively
  hyperbolic groups}, J. Reine Angew. Math. \textbf{612} (2007), 1--15.
  \MR{2364071}

\bibitem[HP97]{HP}
Sa'ar Hersonsky and Fr\'ed\'eric Paulin, \emph{On the rigidity of discrete
  isometry groups of negatively curved spaces}, Comment. Math. Helv.
  \textbf{72} (1997), no.~3, 349--388. \MR{1476054}

\bibitem[KH95]{KH95}
Anatole Katok and Boris Hasselblatt, \emph{Introduction to the modern theory of
  dynamical systems}, Encyclopedia of Mathematics and its Applications,
  vol.~54, Cambridge University Press, Cambridge, 1995. \MR{1326374}

\bibitem[Kle06]{kleiner}
Bruce Kleiner, \emph{The asymptotic geometry of negatively curved spaces:
  uniformization, geometrization and rigidity}, International {C}ongress of
  {M}athematicians. {V}ol. {II}, Eur. Math. Soc., Z\"{u}rich, 2006,
  pp.~743--768. \MR{2275621}

\bibitem[Laf07]{LF07}
Jean-Fran\c{c}ois Lafont, \emph{Diagram rigidity for geometric amalgamations of
  free groups}, no.~3, 771--780.

\bibitem[Liv71]{livsic}
A.~N. Liv\v{s}ic, \emph{Certain properties of the homology of {$Y$}-systems},
  Mat. Zametki \textbf{10} (1971), 555--564. \MR{293669}

\bibitem[Lov75]{lovasz}
L.~Lov\'asz, \emph{Three short proofs in graph theory}, J. Combinatorial Theory
  Ser. B \textbf{19} (1975), no.~3, 269--271. \MR{396344}

\bibitem[LS98]{ljs}
Jouni Luukkainen and Eero Saksman, \emph{Every complete doubling metric space
  carries a doubling measure}, Proc. Amer. Math. Soc. \textbf{126} (1998),
  no.~2, 531--534. \MR{1443161}

\bibitem[LT03]{LT}
Artur~O. Lopes and Philippe Thieullen, \emph{Sub-actions for {A}nosov
  diffeomorphisms}, no. 287, 2003, Geometric methods in dynamics. II, pp.~xix,
  135--146. \MR{2040005}

\bibitem[Ota90]{otal}
Jean-Pierre Otal, \emph{Le spectre marqu\'e{} des longueurs des surfaces \`a{}
  courbure n\'egative}, Ann. of Math. (2) \textbf{131} (1990), no.~1, 151--162.
  \MR{1038361}

\bibitem[Pol87]{Pol87}
Mark Pollicott, \emph{Symbolic dynamics for {S}male flows}, Amer. J. Math.
  \textbf{109} (1987), no.~1, 183--200. \MR{878205}

\bibitem[PS04]{ps}
Mark Pollicott and Richard Sharp, \emph{Livsic theorems, maximizing measures
  and the stable norm}, Dyn. Syst. \textbf{19} (2004), no.~1, 75--88.
  \MR{2038273}

\bibitem[Sig72]{sigmundflow}
Karl Sigmund, \emph{On the space of invariant measures for hyperbolic flows},
  Amer. J. Math. \textbf{94} (1972), 31--37. \MR{302866}

\bibitem[Sig74]{sigmund}
\bysame, \emph{On dynamical systems with the specification property}, Trans.
  Amer. Math. Soc. \textbf{190} (1974), 285--299. \MR{352411}

\bibitem[{Wu}23]{wu}
Yandi {Wu}, \emph{{Marked Length Spectrum Rigidity for Surface Amalgams}},
  arXiv e-prints (2023), arXiv:2310.09968.

\end{thebibliography}
\Addresses

\end{document}